\documentclass[10pt,reqno]{amsart}
\usepackage{times}
\usepackage[T1]{fontenc}
\usepackage{amssymb}
\usepackage{amsmath,amsthm}
\usepackage{amsfonts}
\usepackage{leftidx}
\usepackage{color}
\usepackage{mathrsfs}
\usepackage{enumerate}	
\usepackage{abstract}
\usepackage{stmaryrd}
\usepackage{ulem}
\usepackage{soul}
\usepackage{cancel}
\usepackage{graphicx}
\usepackage{appendix}
 \usepackage{hyperref}
 \usepackage{graphicx} 
\usepackage{float} 
\usepackage{subfigure} 
\usepackage[top=0.8in, bottom=0.7in, left=0.8in, right=0.8in]{geometry}

\def\be{\begin{equation}}
\def\ee{\end{equation}}

\newtheorem{theorem}{Theorem}[section]
 \newtheorem{corollary}{Corollary}[section]

\newtheorem{lemma}{Lemma}[section]
\newtheorem{proposition}{Proposition}[section]
\theoremstyle{definition}

\theoremstyle{remark}
\newtheorem{remark}{Remark}[section]

\numberwithin{equation}{section}
\begin{document}
 \title[Smoluchowski equation]{Patterns in a Smoluchowski equation}
\author{Xingyu Li and Arghir Zarnescu}
\address{(Xingyu LI)\newline Universit\`{a} degli Studi di Trieste, Dipartimento di Matematica, Informatica e Geoscienze, Trieste, Italy
}
\email{xingyuli92@gmail.com
}
\address{(Arghir Zarnescu)\newline BCAM, Basque Center for Applied Mathematics, Mazarredo 14, E48009 Bilbao, Bizkaia, Spain\newline 
IKERBASQUE, Basque Foundation for Science, Maria Diaz de Haro 3, 48013, Bilbao, Bizkaia, Spain\newline
Simion Stoilow Institute of Mathematics of the Romanian Academy,
P.O. Box 1-764, RO-014700 Bucharest, Romania
}
\email{azarnescu@bcamath.org 
}

\thanks{}

\maketitle

\begin{abstract}
We analyze the dynamics of concentrated polymer solutions
modeled by  a 2D Smoluchowski equation. We describe the long time behavior of
the polymer suspensions in a fluid.
\par When the flow influence is neglected the equation has a gradient structure.
The presence of a simple flow introduces significant structural
changes in the dynamics.  We study the case of an externally imposed
flow with homogeneous gradient. We show that the equation is still
dissipative but new phenomena appear. The dynamics depend on both
the concentration intensity and the structure of the flow. In
certain {\it limit cases} the equation has a gradient structure, in
an appropriate reference frame, and the solutions evolve to either a
steady state or a tumbling wave. For small perturbations of the
gradient structure we show that  for small concentrations the solutions evolve in
the long time limit to a steady state. However for  high concentrations there is a rigidity phenomenon for the tumbling wave.

\end{abstract}

\section{Introduction}

\par In this paper we study qualitative properties of  a Smoluchowski  equation
 describing the dynamics of non-Newtonian complex fluids containing
liquid crystalline polymers, in a concentrated regime.
\par The model we use was introduced by M. Doi in {\cite{DOI81}} (see also \cite{DOED86}). It
identifies  the polymers with inflexible rods whose thickness is
much smaller than their length.  We  study here the case in which
the fluid is two dimensional. This represents a simplification that
 preserves many of the qualitative features of the physical three
dimensional phenomenon (\cite{LARS90}). This model as well as its
three dimensional analogue has attracted much interest (see \cite{CFTZ},\cite{CKT04}, \cite {CVT05},
\cite{CV05},\cite{[FZQ05a]},\cite{ION16},\cite{[FZQ05b]},\cite{[FZQ04a]},\cite{[FZQ04b]}
 ,\cite{[JOLB04]},\cite{[OTTA]},\cite{[QSZ05]},\cite{AZ06}). A number of further developments are documented in  \cite{BALL21},\cite{BRE}, \cite{DEG20}, \cite{FRA}, \cite{FRA2}, \cite{NIK19},\cite{MAK},\cite{OS}, \cite{OS1} \cite{REZ}, \cite{SI}, \cite{BEZ15}, \cite{BEZ19}, \cite{BEZ20}, \cite{BEZ21}, \cite{LUP17}, \cite{WANG21}.
 
Furthermore, investigations into the Smoluchowski Equation coupled with Navier-Stokes equations have been pursued in \cite{CMA} and \cite{HUANG}, while related studies on coagulation-diffusion models are available in \cite{BRE1} and \cite{BRE2}. Progress regarding more general polymeric fluids is outlined in \cite{JIANG}, \cite{LA1}, \cite{LA2}, and \cite{MEN}.

\par The local probability measure associated with the polymers is
 of the form $f(t,x,\theta)d\theta$. Here $t$ is the time
coordinate, $x\in\mathbb{R}^2$ denotes  the spatial coordinate and
$\theta\in [0,2\pi]$ is a direction on the unit circle. The measure
$fd\theta$ represents the time-dependent probability that a rod with
center of mass at $x$ has an axis of direction $\theta$ in the area
element $d\theta$. The equation we study, in non-dimensional form,
is:
 \begin{equation}
\partial_t f+ u\cdot \nabla_x f+\partial_\theta(V f+
\partial_\theta(\mathcal{K}f)\cdot f)=\partial_{\theta\theta}f,
\label{eq:f}
\end{equation} where $u$ is the
velocity of the underlying fluid. For $\nabla_x
u=(u^i_j)_{i,j=1,2}$,

\begin{equation}
V(x,\theta)=-u^1_1(x)\cos\theta\sin\theta-u^1_2(x)\sin^2\theta+u^2_1(x)\cos^2\theta+u^2_2(x)
 \sin\theta\cos\theta\label{V}
\end{equation} denotes the projection of  $\nabla_x u\cdot(\cos(\theta),\sin(\theta))^t$ on the
tangent space in $(\cos(\theta),\sin(\theta))^t$. This term
describes the way the fluid influences the evolution of $f$.
\par The term $\mathcal{K}f$ -- the excluded volume potential, which
accounts for the interaction between different rods -- is given by

\begin{equation}
\mathcal{K}f(\theta)=b\int_{S^1}k(\theta-\theta')f(\theta')d\theta',
\label{kf}
\end{equation} where $k$ is a smooth function and $b$ is a
non-dimensional parameter measuring the concentration of the
polymers  in the fluid. Moreover,
\begin{equation}k(\theta)=k(-\theta).
\label{k} \end{equation} \par  In many instances we will restrict
ourselves to using $k(\theta)=\frac{1}{2}\cos(2\theta)$ in which case
$\mathcal{K}f$ is the so called Maier-Saupe potential. This
potential  has been frequently used in the literature (\cite{CKT04}, \cite
{CVT05}, \cite{CV05}).
\par The   fluid  velocity $u$ obeys the Stokes or the Navier-Stokes
equations forced by an appropriate average of $f$ (\cite{DOED86}).

 \par The rich dynamical behavior of the system poses  significant
 numerical and analytical challenges (\cite{MACR95},
 \cite{MAMA89}). We consider two levels of
 complexity:

 \par On a first level, we neglect the influence of the fluid. This scenario was analyzed in previous works (\cite{CKT04}, \cite{CVT05}), where it was shown that the system possesses a gradient structure and is dissipative: solutions starting from an arbitrary initial data reach in the long-time limit a  fixed ball. We refine this analysis by demonstrating that the $\omega$-limit set of any solution consists entirely of steady states. Furthermore, under additional symmetry constraints, the $\omega$-limit set reduces to a single steady state.

\par When fluid flow is introduced, the dynamics grow highly complex, even for flows with simple structures. We focus on an externally imposed flow with a homogeneous gradient, meaning that the matrix $\nabla_x u$ is constant. This setup is physically significant as it captures the local behavior of steady, smooth flows, and is frequently used in rheological studies (see \cite{BCAH87}, Chapter 5). A related case, of small-amplitude shear flow, was examined in \cite{ZHWA06}.
 \par We establish that the equation remains dissipative even in the presence of homogeneous gradient flows. The dynamics depend critically on both the concentration intensity $b$ and  the flow structure. For homogeneous gradient flows, $V$ becomes independent of
$x$ and can be expressed as
$$V=\omega+s\cos(2\theta+\alpha),$$ where $\omega\in\mathbb{R},s\ge 0$ are uniquely determined
and $\alpha$ is defined modulo $2\pi$. Notably, $\omega$ corresponds to the vorticity of the flow, calculated as $\omega=\frac{u^2_1-u^1_2}{2}$. Without loss of generality, we set $\alpha=0$ (if $\alpha\ne 0$, the transformation
$\tilde f(t,\theta)=f(t,\theta-\alpha/2)$ recasts the equation with 
$\cos(2\theta)$ replacing $\cos(2\theta+\alpha)$). Here, $\omega$ and $s$ govern the flows rotational and elongational characteristics, respectively.

\par  In \textit{limit cases}, the system retains a gradient structure: when $\omega=0$, this structure persists in the initial reference frame, while for $s=0$, it holds in a rotating frame. These parameters play distinct roles: $s$ primarily influences the long-time shape of solutions, whereas 
$\omega$ governs rotational behavior.

\par For general flows with both $\omega\ne 0$ and $s\ne 0$, no explicit gradient structure exists, though residual features of gradient dynamics remain. When $\omega=0$ is arbitrary and $s$ and $b$ are  small solutions asymptotically approach steady states.

\par Additionally, in the Maier-Saupe potential case ($\omega=0$), stable solutions attract others with an exponential convergence rate dictated by the spectral gap of the linearized problem. This gap must be computed in a norm adapted to the nonlocal term, a method employed in prior studies on diffusion equations (\cite{CD}, \cite{LI19}).

 However, for $\omega>0$ and small $s$  the time-periodic solution of period proportional to the difference  between $\omega$ and $s$ disappears in a neighbourhood of the periodic solution that exists at $s=0$.  It should be pointed out that this however does not rule out the existence of time-periodic solutions  even for small $s\not=0$, but these would either have different time-periods or be further away from the known time-periodic solution (the one at $s=0$).  

\par  Several open questions remain. For small $s$, our analysis applies to extremes of $b$ (either small or large)  and a fixed-time period but not to intermediate values, precisely the regime where transitions between isotropic and nematic phases occur in flow-free systems. The case of large 
$s$ remains particularly challenging.
\par Beyond specific regimes, broader questions arise about the generality of these dynamics. For 1D \textit{local} nonlinear parabolic equations, bounded trajectories typically approach steady states or periodic solutions (\cite{FIMA89}, \cite{MATA86}). Nonlocal terms, however, may permit richer dynamics (\cite{FIPO90}), leaving uncertainty about whether analogous results (e.g., a Poincar\'e-Bendixson theorem) hold here. Further insights into nematic polymer dynamics can be found in \cite{ZHWA06}, \cite{ZHWA11}.

 \par The paper is organized as follows: in the next section we
 consider a general nonlinear Fokker-Planck equation and prove its
 dissipativity. We  recognize that both levels of
 complexity can be put into this general form and this  gives us
 the dissipativity at both levels. Moreover, when the nonlinear
 Fokker-Planck equation has a gradient structure we  analyze its
 $\omega$-limit sets. We  obtain thus information about the case
 when the flow is neglected or the flow is present but irrotational.
 \par In the last section  we  consider the case when the flow is
 present and  is externally imposed, of constant gradient.
  If $s$  and $b$ are small we prove the evolution to a steady
 state. The proof also shows the existence of a unique steady state by
 using arguments from the dynamical systems theory.
 \par Finally,  in the last section we  prove our rigidity theorem concerning  the non-existence
 of  certain time periodic solutions in the appropriate moving frame for
 small $s$ and large $b$. The proof involves the use of  a
 non-standard form of the implicit function theorem in Banach spaces.

\bigskip
 \section{The dissipativity}
 \label{sec2}
 \subsection{ The general case}

 \par Consider the nonlinear Fokker-Planck
 equation:

 \begin{equation}
 \partial_t f+\partial_\theta(fV+f\partial_\theta\mathcal{K}f)=\partial_{\theta\theta} f,
 \label{eq:fgen}
 \end{equation} where $f$ depends only on time and $\theta\in [0,2\pi]$ with $V$  a potential
  which depends on $\theta\in [0,2\pi]$ and may depend on time i.e. $V=V(t,\theta)$. We assume that $f(t,\cdot)$
   is periodic $f(t,0)=f(t,2\pi),\forall t\ge 0$ and since $f$ represents a probability distribution function we assume that
   for all $t\ge 0$ it is positive  and of mean one. We also assume that $V$ is smooth, periodic in the $\theta$ variable  and
   $\mathcal{K}f$ is defined as in (\ref{kf}).

 \par We prove that any solution  of the equation (\ref{eq:fgen}) which starts from a
 nonnegative, mean one, smooth initial data
  will eventually enter a fixed ball in $H^1$.
 Let us remark that for the case when $V=0$, $\mathcal{K}f$ is the  Maier-Saupe potential
 and even initial data the existence of a global attractor in any Sobolev norm,
 was proved in \cite{CVT05}.

 \smallskip
 \begin{theorem} (The dissipativity) Assume that there exist
 $M,N$ such that
 \begin{displaymath}
 \|\partial_\theta V(t,\theta)\|_{L^\infty(S^1)}\le M,\forall
 t\ge 0\,\, |\partial_{\theta\theta}k|\le N,\,\forall
 \theta.
 \end{displaymath}

 \par Let $f$ be a solution of (\ref{eq:fgen}) starting from a smooth
 nonnegative initial data, with $\|f(0)\|_{L^1}=1$. We have

 \begin{equation}
 \|f(t)\|_{L^2}^2\le \|f(0)\|_{L^2}^2 e^{-t/2}+
 \bar C+\frac{1}{\pi},
 \label{l2diss}
 \end{equation}
 \begin{equation}
 \|\partial_\theta f(t)\|_{L^2}^2\le \left(\|\partial_\theta f(0)\|_{L^2}^2+\|f(0)\|_{L^2}^2\cdot 6(M+bN)\right)e^{-t/4}+(M+bN)\left(8\bar
C+\frac{6}{\pi}\right), \label{h1diss}
 \end{equation}with $$\bar C=4\left(\frac{M+bN}{2}(\frac{1}{2\pi}+C2^{1/3})+
 (\frac{MC}{2^{2/3}})^2+(\frac{bNC}{2^{2/3}})^2\right),$$  where $C$ is a constant that
 appears in the Gagliardo-Nirenberg inequality (\ref{gani}).
 \label{dissipativity}
 \end{theorem}
\begin{proof}The existence of solutions is obtained by standard arguments (see also
 \cite{CKT04}). The positivity of the initial data is
 preserved, by the maximum principle, and since
 $\frac{d}{dt}\int_{S^1}f(t,\theta)d\theta=0$ we have that
 $\|f(t,\cdot)\|_{L^1}=1,\forall t\ge 0$.

\par Let us recall  the classical Gagliardo-Nirenberg inequalities
(\cite{NIRE59}) which state that for $1\le q,r\le \infty$ and $j,m$
integers satisfying $0\le j<m$ and for $f$ an appropriately smooth,
 periodic, mean zero function
\begin{equation} \|D^j f\|_{L^p}\le C\|D^m f\|^a_{L^r}
\|f\|_{L^q}^{1-a},\label{gani}
\end{equation}
with
\begin{displaymath} \frac{1}{p}=\frac{j}{d}+a(\frac{1}{r}-\frac{m}{d})+
\frac{1-a}{q}, \end{displaymath} where $d$ is the spatial dimension.

\par Taking $d=1,j=0,p=\infty,m=1,r=2,q=1$ we obtain

\begin{displaymath}
\|f-\bar f\|_{L^\infty}\le C \|\partial_\theta
f\|_{L^2}^{2/3}\|f-\bar f\|^{1/3}_{L^1}, \end{displaymath} where we
denoted by $\bar f$ the average of $f$. Using the fact that the
$L^1$ norm of $f$ is $1$ and $\bar f=\frac{1}{2\pi}$  the last
inequality implies
\begin{equation}
\|f\|_{L^\infty}\le \frac{1}{2\pi}+C 2^{1/3}(1+\|\partial_\theta
f\|_{L^2}).\label{rel:finfty}
\end{equation}
\par Also,  using Poincar\'{e}'s inequality and
$\bar f=\frac{1}{2\pi}$ we have
\begin{displaymath}
\|f\|_{L^2}\le \|\bar f\|_{L^2}+\|f-\bar f\|_{L^2}\le
\frac{1}{\sqrt{2\pi}}+\|\partial_\theta f\|_{L^2},
\end{displaymath} and then
\begin{equation}
\|f\|_{L^2}^2\le \frac{1}{\pi}+2\|\partial_\theta f\|_{L^2}^2.
\label{rel:fl2a}
\end{equation}
\par Multiplying (\ref{eq:fgen}) by $f$, integrating over $S^1$ and by parts we obtain
\begin{equation}
\frac{1}{2}\frac{d}{dt}\int_{S^1} f^2+\frac{1}{2}\int_{S^1}
f^2\partial_\theta
V+\frac{1}{2}\int_{S^1}(\partial_{\theta\theta}\mathcal{K}f)
f^2=-\int_{S^1} (\partial_\theta f)^2. \label{balance1}
\end{equation}

\par Using the hypothesis, the bound (\ref{rel:finfty})
and $\|f\|_{L^1}=1$ we get

\begin{equation*}
\begin{aligned}
\frac{1}{2}|\int_{S^1} \partial_\theta V f^2d\theta|&\le \frac{1}{2}
\|\partial_\theta V\|_{L^\infty}\|f\|_{L^\infty}
\|f\|_{L^1}\le \frac{M}{2}\left(\frac{1}{2\pi}+C2^{1/3}(1+\|\partial_\theta
f\|_{L^2})\right)\\&\le \frac{M}{2}(\frac{1}{2\pi}+C2^{1/3})+
(\frac{MC2^{1/3}}{2})^2+\frac{1}{4} \|\partial_\theta f\|_{L^2}^2
\nonumber.
\end{aligned}
\end{equation*} \par Similarly

\begin{displaymath}
\frac{1}{2}|\int_{S^1} (\partial_{\theta\theta}\mathcal{K}f)f^2
d\theta|\le \frac{bN}{2}(\frac{1}{2\pi}+C2^{1/3})+
(\frac{bNC2^{1/3}}{2})^2+\frac{1}{4} \|\partial_\theta f\|_{L^2}^2.
\end{displaymath}

\par Using the last two bounds in (\ref{balance1}) we obtain

\begin{equation}
\frac{1}{2}\frac{d}{dt}\|f\|_{L^2(S^1)}^2\le \underbrace{
 \frac{M+bN}{2}(\frac{1}{2\pi}+C2^{1/3})+
 (\frac{MC}{2^{2/3}})^2+(\frac{bNC}{2^{2/3}})^2}_{\bar
 C/4}-\frac{1}{2}\|\partial_\theta f\|_{L^2}^2.
 \label{eq:diss}
 \end{equation} So then using (\ref{rel:fl2a}), multiplying by $2e^{t/2}$,
integrating on $[0,t]$ and then multiplying by $e^{-t/2}$ we obtain
(\ref{l2diss}).

\par On the other hand (\ref{eq:diss}) can be rewritten as
 \begin{displaymath}\frac{1}{2}\frac{d}{dt}\|f\|_{L^2}^2+\frac{1}{4}\|\partial_\theta
 f\|_{L^2}^2\le \frac{\bar C}{4}-\frac{1}{4}\|\partial_\theta f\|_{L^2}^2.
 \end{displaymath}

 \par Using  (\ref{rel:fl2a}) on the right hand side, multiplying by $2e^{t/4}$,
 integrating on $[0,t]$
 and multiplying by $e^{-t/4}$ we get

\begin{eqnarray}
\|f(t)\|_{L^2}^2+\frac{1}{2}\int_0^t \|\partial_\theta
f(s)\|_{L^2}^2 e^{(s-t)/4}ds\le2\bar
C+\frac{1}{\pi}+\|f(0)\|_{L^2}^2e^{-t/4}, \label{P}
\end{eqnarray} which gives us an apriori bound on a time integral
involving $\|\partial_\theta f\|_{L^2}^2$.

\par In order to obtain the dissipativity of $\|\partial_\theta f\|_{L^2}$
let us differentiate (\ref{eq:fgen}) with respect to $\theta$.
Denoting $\partial_\theta f=F$ we obtain an equation for $F$:

\begin{equation}
\partial_t F+\partial_\theta\left[FV+f\partial_\theta V
+F\partial_\theta(\mathcal{K}f)+f\partial_{\theta\theta}
(\mathcal{K}f)\right]=\partial_{\theta\theta}F. \label{F}
\end{equation}

\par Multiplying by $F$, integrating on $S^1$ and by parts we
have

\begin{eqnarray}
\frac{1}{2}\frac{d}{dt}\int_{S^1}F^2+\frac{1}{2}\int_{S^1}
\partial_\theta V F^2-\int_{S^1}\partial_\theta V
f\partial_\theta F\nonumber
+\frac{1}{2}\int_{S^1}\partial_{\theta\theta}(\mathcal{K}f)F^2-
\int_{S^1}f\partial_{\theta\theta}(\mathcal{K}f)\partial_\theta
F=-\int_{S^1} (\partial_\theta F)^2\nonumber.
\end{eqnarray}
\par Hence

\begin{equation*}
\begin{aligned}
\frac{1}{2}\frac{d}{dt}\|F\|_{L^2}^2&=\int_{S^1}(\partial_\theta
V+\partial_{\theta\theta}(\mathcal{K}f))f\partial_\theta F
-\frac{1}{2}\int_{S^1}(\partial_\theta V+\partial_{\theta\theta}(\mathcal{K}f))F^2-\int_{S^1}(\partial_\theta F)^2\nonumber\\
&\le\frac{1}{2} (\|\partial_\theta
V\|_{L^\infty}+b\|\partial_{\theta\theta}k\|_{L^\infty})\int_{S^1}f^2+\frac{1}{2}\int_{S^1}(\partial_\theta
F)^2 \\&+ \frac{1}{2}(\|\partial_\theta
V\|_{L^\infty}+b\|\partial_{\theta\theta}k\|_{L^\infty})\int_{S^1}F^2-\int_{S^1}
(\partial_\theta F)^2\nonumber\\
&\le \frac{(M+bN)}{2}(P+\|F\|_{L^2(S^1)}^2)-\frac{\|\partial_\theta
F\|_{L^2(S^1)}^2}{2}\nonumber,
\end{aligned}
\end{equation*}
 where for the last inequality we denoted $P=\|f(0)\|_{L^2}^2 e^{-t/2}+\bar
C+\frac{1}{\pi}$ and used (\ref{l2diss}) which we have just proved.
\par Using the fact that $F=f_\theta$ is mean zero and Poincar\'{e}'s
inequality the last inequality implies

\begin{displaymath}
\frac{d}{dt}\|F\|_{L^2}^2+\frac{\|F\|_{L^2}^2}{4}\le
(M+bN)(P+\|F\|_{L^2}^2).
\end{displaymath}

\par Multiplying by $e^{t/4}$, integrating on $[0,t]$ and
multiplying by $e^{-t/4}$ we have

\begin{equation*}
\begin{aligned}
\|F(t)\|_{L^2}^2&\le
\left(\|F(0)\|_{L^2}^2+4(M+bN)\|f(0)\|_{L^2}\right)e^{-t/4}\nonumber
+(M+bN)\left(4(\bar C+\frac{1}{\pi})+\int_0^t
\|F(s)\|_{L^2}^2 e^{(s-t)/4}ds\right)\\&
\le (\|F(0)\|_{L^2}^2+6(M+bN)\|f(0)\|_{L^2}^2)e^{-t/4}+(M+bN)(8\bar
C+\frac{6}{\pi}),\end{aligned}
\end{equation*} where for the last inequality we used
(\ref{P}).\end{proof}
\begin{remark}  Repeating the procedures described above for higher order derivatives
 allows one to obtain, inductively, that the equation is dissipative in any Sobolev norm.
\end{remark}
\begin{remark} One can improve the rate of decay into the absorbing
ball at the expense of having an absorbing ball of bigger radius.
\end{remark}

\begin{remark} Note that  equation (\ref{eq:fgen}) generates a compact
nonlinear semigroup. Indeed, consider the mapping $S(t):H^1\cap
\mathcal{C}\cap S_{L^1}(0,1)\to H^1\cap \mathcal{C}\cap
S_{L^1}(0,1)$ which associates to an element $f$ the solution at
time $t$ starting from initial data $f$ (we denoted by $\mathcal{C}$
the cone of nonnegative functions and by $S_{L^1}(0,1)$ the unit
sphere in $L^1(S^1)$). Then, as $\|f(t)\|_{L^1}=1,\forall t\ge 0$
and $\mathcal{K}f$ is the convolution of a smooth kernel $k$ with
$f$, we have that $\mathcal{K}f$ is apriori bounded in any Sobolev
norm. Thus the equation can be essentially treated  as a semilinear
parabolic equation and standard arguments give the existence,
uniqueness and continuous dependence on the initial data which shows
that $S(t)$ is a semigroup. The compactness is a consequence of the
usual smoothing effect of parabolic equations. See for details
\cite{HENRY81}, Ch.3. \label{remark:csemi}
\end{remark}
\subsection{The gradient case}

\par In the following we assume that there exists
a $2\pi$-periodic function $W=W(\theta)$ such that

\begin{equation}
V=\partial_\theta W. \label{gradreq}
\end{equation}

\par We have then that  (\ref{eq:fgen}) becomes an
equation of gradient type with the free energy functional

\begin{equation}
\label{freeenergy0}
\mathcal{E}=\int_{S^1}\log f\cdot
f-\frac{1}{2}\int_{S^1}\mathcal{K}f\cdot f-\int_{S^1}W\cdot f,
\end{equation} (see also \cite{COKT04b}) and the Fisher information

\begin{equation}
\mathcal I:=\frac{d\mathcal{E}}{dt}=-\int_{S^1}|\partial_\theta(\log
f-\mathcal{K}f-W)|^2 f\,d\theta.\label{energydecay}
\end{equation}

\par We  show that the presence of this energy functional is
enough for proving that the $\omega$-limit set of any solution is
made of steady states. We  show that  the
$\omega$-limit  reduces to only one steady state if additional symmetry
 constraints are imposed.

\par We  first  need some properties of the energy functional:
\begin{lemma} Assume that $W\in L^\infty(S^1)$. Then the energy functional
$\mathcal{E}$  is bounded from below along the solutions and it is
locally Lipschitz as a functional from $L^2\cap\mathcal{C}$ into $\mathbb{R}$
(where $\mathcal{C}$ denotes the cone of nonnegative functions).
\label{lemma:energy}
\end{lemma}

\begin{proof}The energy $\mathcal{E}(f)=\int_{S^1} f\log f-\frac{1}{2}\int_{S^1}
\mathcal{K}f\cdot f-\int_{S^1} fW$
 is made of three parts: the
(negative) Boltzmann entropy $\int_{S^1} f\log f$, the nonlinear
potential contribution $\frac{1}{2}\int_{S^1} \mathcal{K}f\cdot f$ and the
linear potential part $\int_{S^1}Wf$.
\par We have that the nonlinear potential contribution part is bounded in
$L^\infty$ thanks to the fact that $f\ge 0$ and $\int_{S^1}f=1$.
Indeed:

\begin{displaymath}
||\int_{S^1} \mathcal{K}f\cdot f||_{L^\infty}\le
||\mathcal{K}f||_{L^\infty}||f||_{L^1}\le ||k||_{L^\infty}
(\int_{S^1} f)^2.
\end{displaymath}

 A similar argument  works for the linear
potential part, using the hypothesis $W\in L^\infty(S^1)$.

\par On the other hand, the function $x\log x$ is bounded from below by $-\frac {1}{e}$
 which combined with the previous
observations gives us the boundedness from below of the energy
$E(f)$.
\par The fact that the entropy is a locally Lipschitz functional in $L^2$ is
a consequence of the  inequality (see also \cite{GAWA05}):

\begin{displaymath}
|x\log x-y\log y|\le C(|x-y|^{\frac
{1}{2}}+|x-y|(x^{\frac{1}{2}}+y^{\frac {1}{2}})).
\end{displaymath}
\par Also, the nonlinear potential part of the energy is locally
Lipschitz in $L^2$ norm:

\begin{eqnarray}
\|\int \mathcal{K}f\cdot f-\mathcal{K}g\cdot g\|_{L^2}\le ||\int
\mathcal{K}f(f-g)+\int \mathcal{K}(f-g)g||_{L^2}\nonumber\le
||\mathcal{K}f||_{L^\infty}||f-g||_{L^2}+||\mathcal{K}(f-g)||_{L^\infty}||g||_{L^2}
\nonumber\\ \le \|\mathcal{K}\|_{L^1\to
L^\infty}\underbrace{\|f\|_{L^1}}_{=1}||f-g||_{L^2}+
\|\mathcal{K}\|_{L^2\to L^\infty}\|f-g\|_{L^2}\|g\|_{L^2}\nonumber,
\end{eqnarray} where we used  the fact that $k$ is smooth and thus $\|\mathcal{K}\|_{L^1\to
L^\infty},\|\mathcal{K}\|_{L^2\to L^\infty}$ are bounded. It is easy
to check that the linear potential part is also locally Lipschitz in
$L^2$.
\par Thus we have that  the entropy and the potential parts of
the energy are locally Lipschitz in $L^2$, which finishes the proof
of the lemma.  \end{proof}

\par We prove now that
 the nonlinear Fokker-Planck equation
evolves to the steady states in the $H^1$ norm.
\begin{lemma} For any nonnegative initial data of the nonlinear Fokker-Planck
equation (\ref{eq:fgen}) satisfying  (\ref{gradreq})  we have that
the $\omega$-limit set of the corresponding trajectory

\begin{displaymath}
\Omega=\{\Psi; f(t_n)\stackrel {H^1}{\rightarrow} \Psi,\mathrm{for\,
some\, sequence} \,(t_n)_{n\in \mathbb{N}}\subset\mathbb{R}_+\}
\end{displaymath} contains only steady states.
\label{lemma:decaytosteady}
\end{lemma}
\begin{proof} The compactness in   $H^1\cap\mathcal{C}$ of
  the semigroup generated by (\ref{eq:fgen}) ( see {\it
  Remark}~\ref{remark:csemi})
together with the fact that all  trajectories decay exponentially
into a fixed ball suffice for having a connected global attractor in
$H^1\cap\mathcal{C}$ (see for instance \cite{HALE88}, p.39, Thm.
3.4.6; note that there one has a semigroup defined on a Banach
space.The fact that the semigroup is invariant with respect to a
cone, as we have here, will not affect the validity of the quoted
result, as one can easily check).

 \par Take an arbitrary nonnegative initial data and consider the omega limit set
 associated to the trajectory starting from this initial data,
 $\Omega$. Observe that all the elements in $\Omega$ have the same
 energy. Indeed, the energy is decreasing along the trajectories and it is
 bounded from below which means that
 there exists a $c\in\mathbb{R}$ so that
 $\lim_{t\to\infty}E(f(t))=c$. We claim  that this implies that
 the elements of $\Omega$ are  the steady states.
 \par  Recall that the  $\omega$-limit set is an invariant set (see \cite{HALE88}, p.36) so any
 trajectory starting from an initial data in $\Omega$ will stay in $\Omega$ and have the same
 energy $c$. For an initial data in $\Omega$ equation (\ref{energydecay}) implies
 $\log f-\mathcal{K}f-W=\mathrm{const},\forall t\ge 0$.
 Indeed, if the right hand side of (\ref{energydecay}) is negative at some time
 $t_0$ then it
  will be negative on an interval around $t_0$ and this would imply that
  $\mathcal{E}(t)<\mathcal{E}(0),\forall t>t_0$, which is a contradiction.  From
  (\ref{eq:fgen}),(\ref{gradreq}) and  $\log f-\mathcal{K}f-V=\mathrm{const},\forall t\ge 0$
  we have  that $f_t(t,\theta)=0,\forall t\ge0, \theta\in
  [0,2\pi]$ i.e. $f(t,\theta)=f(0,\theta),\forall t\ge 0,\theta\in
  [0,2\pi]$. \end{proof}
  \smallskip
\subsection{Convergence and asymptotic behavior to the stationary states with symmetry constraint}
  \par In general we do not know   if $\Omega$
  reduces to only one steady state. However, when we have a certain type of
  additional symmetry  constraint then there are only finitely many steady
   states which can exist in $\Omega$ and only one of
  them  will actually be in $\Omega$ for a given trajectory. We show
  this for the simple case when $k=\frac{1}{2}\cos(2\theta)$ so $\mathcal{K}=\mathcal{K}_{MS}$ and
   $V=0$ hence (\ref{eq:fgen}) reduces to

\begin{equation}
\partial_t f+\partial_\theta(f\partial_\theta\mathcal{K}_{MS}f)
=\partial_{\theta\theta} f, \quad f(0,\cdot)=f_{{int}}.
 \label{eq:fgradms}
 \end{equation}

\par  This equation has been analyzed in \cite{CVT05} where it was
proved that if $f_{int}$ is even, then the
evenness of $f_{int}$ will be preserved by the flow.
Therefore one can restrict oneself to studying solutions which have
this symmetry, i.e. solutions of the form

\begin{displaymath}
f(t,\theta)=\frac{1}{2\pi}+\frac{1}{\pi}\Sigma_{k=1}^\infty y_k(t)
\cos(2k\theta),
\end{displaymath} where
$y_k(t)=\int_0^{2\pi}\cos(2k\theta)f(t,\theta)d\theta$. The normalization of the initial data implies $y_0=1$ and
$|y_k|\le 1$, and the nonlinear interaction potential becomes
$\mathcal{K}_{MS}(\theta,t)=\frac{b}{2}y_1(t)\cos(2\theta)$.

From Theorem 1 of \cite{CV05}, we know that when $0<b\le 4$, (\ref{eq:fgradms})  has only one even steady state \[f_0(\theta):=\frac 1{2\pi}.\]And for $b>4$, from the results of \cite{FAS},
there are three even steady states: an isotropic solution $f_0$ and two nematic solutions \begin{equation}
\label{Bessela}f_{\pm r}(\theta):=\frac{e^{\pm r\cos (2\theta)}}{\int_0^{2\pi}e^{r\cos (2\theta)} d\theta},\end{equation}here $r(b)$ is the unique positive number that satisfies
\begin{equation}
\label{Bessel2}
\frac{rI_0(r)}{I_1(r)}=\frac b2,
\end{equation}
where $I_k(r)$ for $k\in\mathbb N$ is the modified Bessel function of first kind
\begin{equation}
\label{Bessel}
I_k(r):=\sum_{n=0}^\infty\frac 1{\Gamma (n+1)\Gamma (k+n+1)}\left(\frac r2\right)^{2n+k}=\frac 1{2\pi}\int_0^{2\pi}{e^{r\cos(2\theta)}\cos(2k\theta)d\theta},
\end{equation}
see (B.31) of \cite{MAI}. Notice that (\ref{eq:fgradms}) can be
rewritten in terms of Fourier coefficients as an infinite system of
ODE's:
\begin{equation}
\label{ode}
y_0=1,\quad 
y'_k+4k^2 y_k=bky_1(y_{k-1}-y_{k+1}),\quad k=1,2,\dots
\end{equation}
For $k=1$ we have
\begin{equation}
\label{equationy1}
y_1'=y_1(-4+b-by_2),
\end{equation} which implies that if the $y_1(0)=0$ then $y_1$
will be 0 for all times, which means $f(t,\theta)$ converges to $f_0$. On the other hand, $y_1(t)$ is always positive if $y_1(0)>0$. If $f(t,\theta)$ converges to $f_0$ as $t\to\infty$, then $y_2(t)$ converges to 0, which means there exists $t_0\ge 0$, such that $|y_2(t)|\le \frac {b-4}{2b}$. from \eqref{equationy1}, we obtain that for $t\ge t_0$, $y_1'\ge \frac{b-4}2y_1$, and by Gr{\"o}nwall inequality $y_1$ goes to infinity as $t\to\infty$, a contradiction. So $f(t,\theta)$  converges to $f_r$ as $t\to\infty$.  The eveness of initial data $f_{int}$, which is propagated by the flow,
  leads thus to only three possible elements in the $\omega$-limit set $\Omega$. But
$\Omega$ must be a connected set (see \cite{HALE88}) thus it will
necessarily consist of only one element.

Furthermore, we continue by providing a result on the asymptotic behavior of $f(t,\theta)$, namely concerning the convergence towards the stationary states. We only consider the case $r\ge 0$ since $r<0$ is similar.

\begin{theorem}
\label{Thm2}
 Consider the equation (\ref{eq:fgradms}) with initial data $f_{{int}}$.\\
 (1) Recall that $f_0=\frac 1{2\pi}$. For any $b>0$, if $\int_0^{2\pi} f_{{int}}\cos (2\theta)d\theta=0$, then for any $t>0$, we have 
\begin{equation}
\label{asym1}
\int_0^{2\pi}{(f(t,\theta)-f_0)^2 d\theta
}\le (f_{{int}}(\theta)-f_0)^2e^{-32t}.
\end{equation}

(2) If $b\in (0,4)\cup (4,\infty)$ and $\int_0^{2\pi} f_{{int}}\cos (2\theta)d\theta>0$, then there exist positive constants $C, \lambda$, such that for any $t>0$,
\begin{equation}
\label{asym3}
\int_0^{2\pi}{\frac{(f(t,\theta)-f_r)^2}{f_r} d\theta
}\le Ce^{-\lambda t},
\end{equation}
where $f_r=f_0$ for $0<b<4$, and $f_r$ is defined as \eqref{Bessela},  and $r$ is given by \eqref{Bessel2} for $b>4$.
\end{theorem}
\begin{proof} (Proof of (1) of Theorem \ref{Thm2})
If $y_1$ is always zero, we directly deduce from \eqref{ode} that $y_k'=-4k^2y_k$ for $k\ge 2$, so\[ \frac{d}{dt}\sum_{k=2}^\infty y_k^2=2\sum_{k=2}^\infty y_ky_k'=-8\sum_{k=2}^\infty k^2 y_k^2\le -32\sum_{k=2}^\infty y_k^2,\] 
and \eqref{asym1} is proved by using Gr{\"o}nwall inequality. 
\end{proof}
The rest of this subsection is devoted to proving \eqref{asym3}. From now on, we suppose that $\int_0^{2\pi} f_{{int}}\cos (2\theta)d\theta>0$. 
\subsubsection{Case 1: $b>4$.}For convenience, we define \begin{equation}\label{zk}z_k:=\int_0^{2\pi}f_r\cos(2k\theta)d\theta=\frac{I_k(r)}{I_0(r)}.\end{equation}
From \eqref{ode},  by letting $t\to\infty$, we obtain $4kz_k=bz_1(z_{k-1}-z_{k+1})$ for $k=1,2...$ (this can be also directly deduced from the properties of Bessel function). Let $k=1$, from the definition of $r$ and \eqref{equationy1},we deduce that 
\begin{equation}
\label{ode1}z_1=\frac {2r}b, \quad  z_2=1-\frac 4b.\end{equation}

Recall the free energy $\mathcal E$ and Fisher information $\mathcal I$ as defined in \eqref{freeenergy0} and \eqref{energydecay}. Here they have the forms
\begin{equation}
\label{energy1}
\mathcal E(f):=\int_0^{2\pi}f\log fd\theta-\frac 12\int_0^{2\pi}{f\mathcal Kfd\theta}=\int_0^{2\pi}f\log fd\theta-\frac b4\left(\int_0^{2\pi}{f\cos(2\theta)d\theta}\right)^2,
\end{equation}
and
\begin{equation}
\label{Fisher1}
\mathcal I[f]:=-\frac {d}{dt}E(f)=\int_0^{2\pi}|\partial_\theta (\log f-\mathcal Kf)|^2f d\theta.
\end{equation}
Recall that $\mathcal{E}$ is bounded from below, and the stationary solution $f_r$ is the minimizer. We need to study the quadratic forms associated with the expansion of $\mathcal E, \mathcal I$ around $f_r$. For a smooth perturbation $g$ of $f_r$ such that $\int_0^{2\pi} gf_rd\theta=0$, we define
\begin{equation}
\label{equationQextra1}
\begin{aligned}
Q_1(g):&=\lim_{\varepsilon\to 0}\frac 2{\varepsilon^2}(\mathcal E(f_r(1+\varepsilon g))-\mathcal E(f_r))\\
&=\lim_{\varepsilon\to 0}\frac 2{\varepsilon^2}(\int_0^{2\pi} f_r(1+\varepsilon g)\log (f_r(1+\varepsilon g))-f_r\log f_r d\theta)\\
&-\frac b2\lim_{\varepsilon\to 0}\frac 1{\varepsilon^2}\left[\left(\int_0^{2\pi}f_r(1+\varepsilon g)\cos(2\theta)d\theta\right)^2-\left(\int_0^{2\pi}f_r\cos(2\theta)d\theta\right)^2\right].
\end{aligned}
\end{equation}
Note that 
\[
 f_r(1+\varepsilon g)\log (f_r(1+\varepsilon g))-f_r\log f_r =f_r(1+\varepsilon g)\log(1+\varepsilon g)+\varepsilon f_rg\log f_r,
\]
and recall that $\log(1+\varepsilon g)\sim \varepsilon g-\frac 12 \varepsilon^2g^2$,  $\log f_r=r\cos (2\theta)+const$ and $\int_0^{2\pi}gf_rd\theta=0$, we have
\begin{equation}
\label{equationQextra2}
\begin{aligned}
Q_1(g)&=\lim_{\varepsilon\to 0}\frac 2{\varepsilon ^2}\int_0^{2\pi}\frac {\varepsilon^2}2 f_rg^2+\varepsilon rf_rg\cos(2\theta) d\theta\\
&-\frac b2\lim_{\varepsilon\to 0}\frac 1{\varepsilon^2}(2\varepsilon\int_0^{2\pi}f_r\cos(2\theta)d\theta\int_0^{2\pi}gf_r\cos(2\theta)d\theta+\varepsilon^2(\int_0^{2\pi}gf_r\cos(2\theta)d\theta)^2)\\
&=\int_0^{2\pi}g^2f_rd\theta-\frac b2\left(\int_0^{2\pi}gf_r\cos(2\theta)d\theta\right)^2\\
&+\lim_{\varepsilon\to 0}\frac 1{\varepsilon}\underbrace{\left[2r\int_0^{2\pi}gf_r\cos(2\theta)d\theta-b\int_0^{2\pi}f_r\cos(2\theta)d\theta\cdot \int_0^{2\pi}gf_r\cos(2\theta)d\theta\right]}_{=0, \,\,\text{because} \,\,\int_0^{2\pi}f_r\cos(2\theta)d\theta=\frac {2r}b}\\
&=\int_0^{2\pi}g^2f_rd\theta-\frac b2\left(\int_0^{2\pi}gf_r\cos(2\theta)d\theta\right)^2.\\
\end{aligned}
\end{equation}
Similarly, we have
\begin{equation}
\label{equationQextra3}
\begin{aligned}
Q_2(g):&=\lim_{\varepsilon\to 0}\frac 1{\varepsilon^2}\mathcal I(f_r(1+\varepsilon g))\\&=\lim_{\varepsilon\to 0}\frac 1{\varepsilon^2}\int_0^{2\pi}|\partial_\theta (\log [f_r(1+\varepsilon g)]-\mathcal K [f_r(1+\varepsilon g)])|^2f_r(1+\varepsilon g) d\theta\\
&=\lim_{\varepsilon\to 0}\frac 1{\varepsilon^2}\int_0^{2\pi}|\partial_\theta \log f_r-\partial_\theta\mathcal Kf_r+\partial_\theta \log(1+\varepsilon g)-\varepsilon\partial_\theta \mathcal K (f_rg)|^2f_rd\theta\\
&=\int_0^{2\pi}|\partial_\theta(g-\mathcal K(f_rg))|^2 f_rd\theta.\,\, (\text{here we used that}\,\,\partial_\theta \log f_r-\partial_\theta\mathcal Kf_r=0 \,\,\text{and}\,\, \log(1+\varepsilon x)\sim \varepsilon x \,\,\text{as} \,\,\varepsilon \to 0)
\end{aligned}
\end{equation}

The next lemma provides coercivity results for $Q_1$ and $Q_2$ and is useful for studying the linear stability around $f_r$.
\begin{lemma}
\label{Poincare}
For any function $g$ that satisfies $\int_0^{2\pi} gf_r d\theta=0$,\\
(1) There exists $\eta(b)>0$, such that $Q_1(g)\ge \eta(b)\int_0^{2\pi}g^2f_rd\theta$.\\
(2) There exists $\zeta(b)>0$, such that $Q_2(g)\ge \zeta(b) Q_1(g)$.
\end{lemma}
\begin{proof}
Noticing that $\int_0^{2\pi} gf_rz_1 d\theta=0$, we get from Cauchy-Schwartz inequality,
\[
\left(\int_0^{2\pi}gf_r\cos(2\theta)d\theta\right)^2\le \int_0^{2\pi}g^2f_rd\theta \cdot\int_0^{2\pi}(\cos(2\theta)-z_1)^2f_rd\theta.
\]
Next, if we could show that $\xi(b):=\int_0^{2\pi}(\cos(2\theta)-z_1)^2f_rd\theta<\frac 2b$,
then from the definition of $Q_1$, we would have
\[
Q_1(g)\ge (1-\frac b2\cdot\xi(b))\int_0^{2\pi}g^2f_rd\theta,
\]
and we could choose $\eta(b)=1-\frac b2\cdot\xi(b)$. Next, recalling \eqref{ode1}, we have
\begin{equation}
\label{xi1}
\xi(b)=\int_0^{2\pi}\cos^2(2\theta)f_rd\theta-z_1^2=\frac 12\int_0^{2\pi}(1+\cos(4\theta))f_rd\theta-z_1^2=\frac 12(1+z_2)-z_1^2=1-\frac 2b-\frac{4r^2}{b^2}.
\end{equation}
So to prove $\xi(b)<\frac 2b$, we only need to show that $1-\frac 4b-\frac{4r^2}{b^2}<0$, which equals to \begin{equation}\label{inbr}b^2-4b<4r^2.\end{equation} Recall from \eqref{Bessel2} that $b=\frac{2rI_0(r)}{I_1(r)}$, so finally it equals to show that for any $r>0$, \begin{equation}
\label{xi2}\frac{I^2_0(r)}{I^2_1(r)}-\frac{2I_0(r)}{rI_1(r)}<1,\end{equation}
which holds provided that we have
\[
\frac{rI_0(r)}{I_1(r)}<1+\sqrt{1+r^2}.
\] The last relationship can be proved by using properties of the Bessel functions (see Proposition 7 of \cite{Yang} for example).

Next, recall that $f_r$ satisfies Poincar\'e inequality \footnote{In general, a sufficient condition for a weight function $w$ satisfies Poincar\'e inequality is that $w$ belongs to Muckenhoupt weights $A_2$, which means that there exists a unversal constant $C>0$, such that for any interval $B$ inside $[0, 2\pi]$, $\frac 1{|B|}\int_B w(x) dx \cdot \frac 1{|B|}\int_B w^{-1}(x) dx \le C$. for $w=f_r$, notice that $\frac{e^{-r}}{2\pi I_0(r)}\le f_r(\theta)\le \frac{e^r}{2\pi I_0(r)}$ for any $\theta\in [0, 2\pi]$, we could choose $C=e^{2r}$. More details can be seen in Chapter 15 of \cite{Heinonen}.}, which means that there exists a constant $p(b)>0$, such that for any function $h$,
\begin{equation}\label{ineq:Poincareweight}
\int_0^{2\pi }|\partial_\theta h|^2f_rd\theta\ge p\left(\int_0^{2\pi}h^2f_rd\theta-\left(\int_0^{2\pi}hf_rd\theta\right)^2\right).
\end{equation}
Recall that $\mathcal K(f_rg)=\frac b2\cos(2\theta)\int_0^{2\pi}{gf_r\cos(2\theta')d\theta'}$. Define $c(g):=\int_0^{2\pi}{gf_r\cos(2\theta)d\theta}$ for simplicity. Then
\begin{equation}
\label{equationQextra4}
\begin{aligned}
\frac 1 p Q_2(g)&\ge \int_0^{2\pi}(g-\mathcal K(f_rg))^2f_rd\theta-\left(\int_0^{2\pi}(g-\mathcal K(f_rg))f_rd\theta\right)^2\\
&= \int_0^{2\pi}(g-\frac b2c(g)\cos(2\theta))^2f_rd\theta-\left(\int_0^{2\pi}\frac b2c(g)\cos(2\theta)f_rd\theta\right)^2(\text{using that} \int_0^{2\pi} gf_rd\theta=0)\\
&=\int_0^{2\pi}g^2f_rd\theta -b\cdot c(g)\int_0^{2\pi}{gf_r\cos(2\theta)}d\theta+\frac{b^2}4c^2(g)\int_0^{2\pi}\cos^2(2\theta)f_rd\theta-\frac{b^2}4c^2(g)\left(\int_0^{2\pi}\cos(2\theta) f_rd\theta\right)^2\\
&=\int_0^{2\pi}g^2f_rd\theta -b\cdot c^2(g)+\frac{b^2}4c^2(g)(1-\frac 2b)-\frac{b^2}4c^2(g)\frac{4r^2}{b^2}\\& (\text{recall that} \int_0^{2\pi}\cos(2\theta)f_rd\theta=\frac {2r}b, \quad  \int_0^{2\pi}\cos^2(2\theta)f_rd\theta=1-\frac 2b)\\
&=\int_0^{2\pi}g^2f_rd\theta+(\frac {b^2}4-\frac 32b-r^2)c^2(g) \\
&=\int_0^{2\pi}g^2f_rd\theta+\frac {b^2}4\left(1-\frac 6b-\frac{4r^2}{b^2}\right)\left(\int_0^{2\pi}gf_r\cos(2\theta)d\theta\right)^2: =\tilde Q_1[g].
\end{aligned}
\end{equation}
Recall that \[\int_0^{2\pi}g^2f_rd\theta\ge (1-\frac 2b-\frac{4r^2}{b^2})\left(\int_0^{2\pi}gf_r\cos(2\theta)d\theta\right)^2= (1-\frac 2b-\frac{4r^2}{b^2})c^2(g).\] So after direct computation, we obtain that $\tilde Q_1[g]\ge \frac{4r^2+4b-b^2}{2b}Q_1[g]$. This means $Q_2(g)\ge\frac{4r^2+4b-b^2}{2b}p Q_1(g)$. \end{proof}
Now we come to prove the large time asymptotic behaviour. First, we introduce a nonlocal scalar product for the linearized evolution operator. For functions $g_1$ and $g_2$ that satisfy $\int_0^{2\pi} g_1f_rd\theta=\int_0^{2\pi} g_2f_rd\theta=0$, define 
\begin{equation}
\label{scalar}
\langle g_1, g_2\rangle:=\int_0^{2\pi}{g_1g_2f_rd\theta}-\frac b2{\int_0^{2\pi}g_1f_r\cos(2\theta)d\theta}\cdot{\int_0^{2\pi}g_2f_r\cos(2\theta)d\theta}.
\end{equation}
Then $\langle g,g\rangle=Q_1(g)$. Next, for the equation (\ref{eq:fgradms}), set $f=f_r(1+g)$. Then (\ref{eq:fgradms}) can be rewritten as
\begin{equation}
\label{scalar2}
\partial_tg=\mathcal Lg-\frac 1{f_r}\partial_\theta(f_rg\partial_\theta\mathcal K(f_rg)),
\end{equation}
where $\mathcal L$ is the linear operator
\begin{equation}
\label{scalar3}
\mathcal Lg:=\frac 1{f_r}\partial_{\theta\theta}(f_rg)-\frac 1{f_r}\partial_\theta(f_rg\partial_\theta\mathcal K(f_r))-\frac 1{f_r}\partial_\theta(f_r\partial_\theta\mathcal K(f_rg))=\frac 1{f_r}\partial_\theta(f_r\partial_\theta(g-\mathcal K(f_rg))).
\end{equation}
We next show that 
\begin{lemma}
$-\langle \mathcal Lg, g\rangle =Q_2(g)$.
\end{lemma}
\begin{proof}Recall that  $\frac b2\int_0^{2\pi}gf_r\cos(2\theta')d\theta' \cdot \cos(2\theta)=\mathcal K(f_rg)$. So we have 
\begin{equation*}
\begin{aligned}
\langle \mathcal Lg, g\rangle &=\int_0^{2\pi}\partial_\theta(f_r\partial_\theta(g-\mathcal K(f_rg)))g d\theta-\frac b2\int_0^{2\pi}{gf_r\cos(2\theta)d\theta}\int_0^{2\pi}\partial_\theta(f_r\partial_\theta(g-\mathcal K(f_rg))) \cos(2\theta) d\theta\\
&=\int_0^{2\pi}\partial_\theta(f_r\partial_\theta(g-\mathcal K(f_rg)))(g-\mathcal K(f_rg))d\theta =-\int_0^{2\pi}|\partial_\theta(g-\mathcal K(f_rg))|^2f_rd\theta=-Q_2(g),
\end{aligned}
\end{equation*}
where we can integrate by parts because all the functions here have the same value on 0 and $2\pi$. \end{proof}
Finally, we prove Theorem \ref{Thm2} by showing \eqref{asym3}.
\begin{proof}(Proof of (2) of Theorem \ref{Thm2})
From the lemmas above, we have
\begin{equation*}
\begin{aligned}
\frac 12 \frac{d}{dt}Q_1(g)&=\frac 12\frac d{dt}\langle g,g\rangle =\langle g, \mathcal Lg\rangle -\langle g, \frac 1{f_r}\partial_\theta(f_rg\partial_\theta\mathcal K(f_rg))\rangle\\
&=-Q_2(g)-\int_0^{2\pi} \partial_\theta(f_rg\partial_\theta\mathcal K(f_rg))(g-\mathcal K(f_rg))d\theta\\
&=-Q_2(g)+\int_0^{2\pi} f_rg\partial_\theta\mathcal K(f_rg)\partial_\theta(g-\mathcal K(f_rg))d\theta\\
&\le -Q_2(g)+\left(\int_0^{2\pi} |\partial_\theta(g-\mathcal K(f_rg))|^2f_rd\theta\right)^\frac 12 \left(\int_0^{2\pi}g^2|\partial_\theta\mathcal K(f_rg)|^2f_rd\theta\right)^\frac 12.
\end{aligned}
\end{equation*}
Note that
\[
|\partial_\theta\mathcal K(f_rg)|=b\left|\sin(2\theta)\int_{0}^{2\pi} gf_r\cos(2\theta)d\theta\right|\lesssim Q_1^\frac 12(g),
\]
so
\[
\left(\int_0^{2\pi}g^2|\partial_\theta\mathcal K(f_rg)|^2f_rd\theta\right)^\frac 12\lesssim ||\partial_\theta\mathcal K(f_rg)||_{L^\infty}\left(\int_0^{2\pi}g^2f_rd\theta\right)^\frac 12\lesssim Q_1(g)^\frac 12\cdot Q_1(g)^\frac 12=Q_1(g).
\]
Hence there exists a constant $C>0$, such that
\begin{equation}
\label{gron1}
\frac 12 \frac{d}{dt}Q_1(g)\le -Q_2(g)+CQ_1(g)Q_2^\frac 12(g).
\end{equation}
Recall that $Q_1(g)\to 0$ as $t\to\infty$ and $Q_2(g)\ge \zeta Q_1(g)$. For any $\alpha>0$, consider the function  $\mathsf f(x)=-x+C\alpha x^\frac 12 (x>0)$, then $\mathsf f$ is decreasing if $-1+\frac 12C\alpha x^{-\frac 12}<0$. Because $Q_1\to 0$ as $t\to\infty$, there exists $T^*>0$, such that $Q_1(g(t))<\frac {4\zeta}{C^2}$ for $t>T^*$. This means that for $t>T^*$, we have
\[
-1+\frac 12CQ_1(g)Q_2(g)^{-\frac 12}\le -1+\frac 12C\zeta^{-\frac 12}Q_1(g)^\frac 12<0.
\]
Come back to \eqref{gron1}, we have for $t> T^*$, \[\frac d{dt}Q_1(g)\le -2\zeta Q_1(g)+2C\zeta^\frac 12 Q_1(g)^\frac 32. \] Thus from Gr{\"o}nwall inequality, there exists a constant $C'>0$, such that for any $t>0$,
\[
Q_1(g)\le C'e^{-2\zeta t},
\]
which finishes the proof of \eqref{asym3}. \end{proof}
\subsubsection{Case 2: $0<b<4$} 
Recall that $f_0=\frac 1{2\pi}$. We define quadratic forms $Q_1, Q_2$ as \eqref{equationQextra1} to \eqref{equationQextra3}. For the perturbation $g$ of $f_0=\frac 1{2\pi}$ such that $\int_0^{2\pi}g d\theta=0$, we define
\[
Q_1(g)=\frac 1{2\pi}\int_0^{2\pi} g^2d\theta-\frac {b}{8\pi^2}\left(\int_0^{2\pi} g\cos(2\theta)d\theta\right)^2,\quad Q_2(g)=\frac 1{2\pi}\int_0^{2\pi}|\partial_\theta(g-\frac 1{2\pi}\mathcal K(g))|^2 d\theta.
\]
Next, we show the coercivity result as in Lemma \ref{Poincare}.  Because
\begin{equation}
\label{use1}
\left(\int_0^{2\pi} g\cos(2\theta) d\theta\right)^2\le \int_0^{2\pi}\cos^2(2\theta)d\theta \cdot \int_0^{2\pi}g^2d\theta=\pi \int_0^{2\pi}g^2d\theta,
\end{equation}
so from $0<b<4$, we obtain that there exists $\eta_1(b)>0$, such that $Q_1(g)\ge \eta_1(b)\int_0^{2\pi} g^2 d\theta$. Next, we still define $c(g)=\frac 1{2\pi}\int_0^{2\pi} g\cos(2\theta) d\theta$ for simplicity. Similarly as \eqref{ineq:Poincareweight}, \eqref{equationQextra4}, there exists $p_1(b)>0$, such that
\[
\frac 1{p_1}Q_2(g)\ge \frac 1{2\pi}\int_0^{2\pi}g^2 d\theta+(\frac {b^2}8-b) c^2(g).
\]
From \eqref{use1}, we have $ \frac 1{2\pi}\int_0^{2\pi}g^2 d\theta\ge 2c^2(g)$. So finally, we have
\[
\frac 1{p_1}Q_2(g)\ge (\frac {b^2}{16}-\frac b2+1)\frac 1{2\pi}\int_0^{2\pi}g^2 d\theta=\frac {(b-4)^2}{32\pi}\int_0^{2\pi}g^2 d\theta.
\]
From the preceding results, we define the scalar product $\langle\cdot, \cdot \rangle$ and the linear operator $\mathcal L$ as introduced earlier. Furthermore, the identities $\langle g,g\rangle= Q_1(g)$ and $\langle \mathcal Lg, g\rangle=-Q_2(g)$  hold, consistent with the established framework. The nonlinear term, being of higher order, is bounded by  $Q_1$ and $Q_2$, which enables us to derive inequality  \eqref{asym3} from Gr{\"o}nwall inequality. For brevity, we omit the detailed proofs here.

\section{Homogeneous gradient flows}

\par In this section we assume that the flow is externally imposed and
of homogeneous gradient, that is:

$$u(x_1,x_2)=(\Sigma_{i=1}^2 u^1_i x_i,\Sigma_{i=1}^2 u^2_i x_i),$$ with
$\{u^i_j\}_{i,j=1,2}$ arbitrary constants.

\par We consider the moving frame transformation

\begin{displaymath}
\bar f(t,x_1,x_2,\theta)=f(t,x_1+t(\Sigma_{i=1}^2 u^1_i x_i),
 x_2+t(\Sigma_{i=1}^2 u^2_i x_i),\theta),
\end{displaymath} and then the equation for $f$, (\ref{eq:f}), becomes an
equation for $\bar f$

\begin{displaymath}
\partial_t \bar f+\partial_\theta[(\omega+s\cos(2\theta+\alpha))\bar f+
\partial_\theta\mathcal{K}\bar f\cdot \bar f]=\partial_{\theta\theta}\bar f,
\end{displaymath} where we used the fact that $V$ as given in (\ref{eq:f}) can be written
$V=\omega+s\cos(2\theta+\alpha)$ with
$\omega=\frac{u^2_1-u^1_2}{2}$  the vorticity of the flow,
$s>0$  uniquely determined and $\alpha$ uniquely determined modulo $2\pi$. The constants
$\omega,s,\alpha$  can be expressed in terms of $u^i_j,i,j=1,2$. One can
easily check  that  $\bar f(t,\theta-\alpha/2)$ satisfies  the
 equation

\begin{equation}
\partial_t \bar f+\partial_\theta[(\omega+s\cos(2\theta))\bar f+
\partial_\theta\mathcal{K}\bar f\cdot \bar f]=\partial_{\theta\theta}\bar f.
\label{eq:fred}
\end{equation}

\par We continue working with the last equation. We start by
observing that this is an equation of the form (\ref{eq:fgen}) and
the assumptions of {\it Theorem}~\ref{dissipativity} are fulfilled.
Thus equation (\ref{eq:fred}) is dissipative.
\par Moreover, we can prove that for arbitrary $\omega$, with small enough $s$ and $b$ the
solutions evolve, in the long time limit, to a steady state. The
strategy of the proof also shows that in the parameter regime given
by assumptions (\ref{condl2}),(\ref{condh1}) below, there exists a
unique steady state solution of (\ref{eq:fred}).

\begin{theorem} Assume that $\omega\in\mathbb{R}$ and $s,b$ are small enough so that

\begin{equation} 1>s+b\sqrt{2\pi}(\bar C+\frac{1}{\pi}+\epsilon)\|\partial_\theta k\|_{L^\infty}+
\frac{b}{2}\|\partial_{\theta\theta}k\|_{L^\infty}, \label{condl2}
\end{equation}

\begin{equation}
\label{condh1}
1>7s+\frac{b}{2}\|\partial^4_\theta
k\|_{L^\infty}+\frac{3}{2}b\|\partial^2_\theta
k\|_{L^\infty}
+b\sqrt{2\pi}\|k\|_{L^\infty}\left(\epsilon+(2s+b\|\partial_{\theta\theta}
k\|_{L^\infty})(8\bar C+\frac{6}{\pi})+\bar
C+\frac{1}{\pi}\right),
\end{equation}
 for some $\epsilon>0$ with $\bar C$ defined in  {\it
Theorem} ~\ref{dissipativity}.
\par Then any solution evolves to the unique steady state as $t\to\infty$.
\end{theorem}
\begin{proof} We prove the statement in two steps. First we
consider the difference between two arbitrary solutions and we show
that after a certain time $t_0$ depending on the size of the initial
data the two solutions will approach each other at an exponential
rate. In the second step we use step $1$ and a contraction argument
to show that in fact any solution will have to evolve to a steady
state.

\smallskip\par{\it Step 1} Consider the difference between two solutions
$f$ and $g$ starting from the initial data $f(0)$ respectively
$g(0)$:

\begin{equation}
\partial_t(f-g)+\omega\partial_\theta (f-g)+\partial_\theta[
s\cos(2\theta)(f-g)+\partial_\theta \mathcal{K}(f-g)\cdot f+
\partial_\theta \mathcal{K}g\cdot
(f-g)]=\partial_{\theta\theta}(f-g).\label{difference}
\end{equation}

\par Multiply by $f-g$, integrate over $S^1$ and by parts:

\begin{eqnarray}
\frac{1}{2}\frac{d}{dt}\int_{S^1} (f-g)^2-s\int_{S^1}
\sin(2\theta)(f-g)^2(\theta)d\theta- \int_{S^1}
\partial_\theta\mathcal {K}(f-g)\cdot f\partial_\theta(f-g)
\nonumber\\+ \frac{1}{2}\int_{S^1}\partial_{\theta\theta}
(\mathcal{K}g)(f-g)^2=-\int_{S^1}\partial_\theta(f-g)^2.
\label{rel:difbalance1}
\end{eqnarray}

\par We have the following bound

\begin{eqnarray}
|\int_{S^1}\partial_\theta\mathcal{K}(f-g)\cdot
f\partial_\theta(f-g)|\le
 \|\partial_\theta \mathcal{K}(f-g)\|_{L^\infty}
 \|f\|_{L^2}\|\partial_\theta(f-g)\|_{L^2}\nonumber\\
 \le\|\partial_\theta\mathcal{K}\|_{L^2\to L^\infty}\|f-g\|_{L^2}
 R_1
 \|\partial_\theta(f-g)\|_{L^2}\le
 b\sqrt{2\pi}\|\partial_\theta k\|_{L^\infty}R_1
 \|\partial_\theta(f-g)\|_{L^2}^2\nonumber,
 \end{eqnarray} where $R_1=\bar C+\frac{1}{\pi}+\epsilon$ (see {\it
 Theorem} ~\ref{dissipativity} for the definition of $\bar C$) and
 the above inequality holds after the time $t_0$ when the solution
 starting from initial data $f(0)$ enters the ball of radius $R_1$
 in $L^2$ (this time $t_0$ depends only on the size of $f(0)$,
 see {\it Theorem}~\ref{dissipativity}).
 For the last inequality we used the fact that $f-g$ is mean
 zero and  Poincar\'{e}'s inequality.

\par Also

\begin{displaymath}
\frac{1}{2}|\int_{S^1}\partial_{\theta\theta}(\mathcal{K}g)(f-g)^2|\le
\frac{1}{2}\|\partial_{\theta\theta}\mathcal{K}\|_{L^1\to
L^\infty}\underbrace{\|g\|_{L^1}}_{=1}\|f-g\|_{L^2}^2\le
\frac{b}{2}\|\partial_{\theta\theta} k\|_{L^\infty}\|f-g\|_{L^2}.
\end{displaymath}

\par Using the above bounds in (\ref{rel:difbalance1}) we obtain,
for $t\ge t_0$,

\begin{eqnarray}
\frac{1}{2}\frac{d}{dt}\|f-g\|_{L^2}^2\le
-\|\partial_\theta(f-g)\|_{L^2}^2+s\|f-g\|_{L^2}^2+
b\sqrt{2\pi}R_1\|\partial_\theta k\|_{L^\infty}
\|\partial_\theta(f-g)\|_{L^2}^2\nonumber+\frac{b}{2}
\|\partial_{\theta\theta} k\|_{L^\infty}\|f-g\|_{L^2}\nonumber,
\end{eqnarray} and by assumption (\ref{condl2}) and Poincar\'{e}'s inequality we have that the
difference between $f$ and $g$ will decay exponentially after time
$t_0$.

\par In order to evaluate the difference in the $H^1$ norm we
take the derivative of the equation (\ref{eq:fred})  with respect to
$\theta$, for two solutions $f$ and $g$. We denote $\partial_\theta
f=F,\partial_\theta g=G$ and then we have an equation for $F-G$:

\begin{eqnarray}
\partial_t(F-G)+\partial_\theta[
(F-G)(\omega+s\cos(2\theta))-2s\sin(2\theta)(f-g)\nonumber
+(F-G)\partial_\theta\mathcal{K}g+F\partial_\theta\mathcal{K}(f-g)\nonumber\\
+(f-g)\partial_{\theta\theta}(\mathcal{K}g)+f\partial_{\theta\theta}\mathcal{K}(f-g)]=\partial_{\theta\theta}(F-G)
\nonumber.
\end{eqnarray}

\par Multiplying by $F-G$, integrating over $S^1$ and by parts we
obtain

\begin{eqnarray}
\frac{1}{2}\frac{d}{dt}\int_{S^1}(F-G)^2-\underbrace{s\int_{S^1}\sin(2\theta)(F-G)^2}_{T_1}
-\underbrace{4s\int_{S^1}\cos(2\theta)(f-g)(F-G)}_{T_2}
-\underbrace{2s\int_{S^1}\sin(2\theta)(F-G)^2}_{T_3}\nonumber\\
+\underbrace{\frac{1}{2}\int_{S^1}\partial_{\theta\theta}(\mathcal{K}g)(F-G)^2}_{T_4}-
\underbrace{\int_{S^1}F\mathcal{K}(F-G)\partial_\theta(F-G)}_{T_5}
+\underbrace{\int_{S^1}(F-G)^2\partial_{\theta\theta}(\mathcal{K}g)}_{T_6}
+\underbrace{\int_{S^1}(f-g)\partial_{\theta}^3(\mathcal{K}g)(F-G)}_{T_7}\nonumber\\
-\underbrace{\int_{S^1}f\partial_\theta\mathcal{K}(F-G)\partial_\theta(F-G)}_{T_8}=-\int_{S^1}\partial_\theta(F-G)^2,
\label{big}
\end{eqnarray} where  we used (\ref{k}) on the last line three lines to commute $\partial_\theta$ and $\mathcal{K}$.

\par We bound each term $T_i,i=1,\dots, 8$ by
$c_i\|F-G\|_{L^2}^2$ or $c_i\|\partial_\theta(F-G)\|_{L^2}^2$ with
some appropriate constants $c_i$.

\begin{equation}
|T_1+T_3+T_4+T_6|\le (s+2s+\frac{3}{2}\cdot
b\|\partial_{\theta\theta}k\|_{L^\infty})\|F-G\|_{L^2}^2.
\label{1346}
\end{equation}

\par We bound $T_2$ and $T_7$ in the same manner. We show only how
to bound $T_2$. Using an integration by parts twice we have

\begin{eqnarray}
|4s\int_{S^1}\cos(2\theta)(f-g)\partial_\theta(f-g)|=
|4s\int_{S^1}\sin(2\theta)(f-g)^2|
\le 4s\|f-g\|_{L^2}^2\le 4s\|F-G\|_{L^2}^2\nonumber,
\end{eqnarray} where we used Poincar\'{e}'s inequality in the last
relation.
\par Then

\begin{equation}
|T_2+T_7|\le (4s+\frac{b}{2}\|\partial^4_\theta k\|_{L^\infty})\|F-G\|_{L^2}^2.
\label{27}
\end{equation}

\par The terms $T_5$ and $T_8$ are treated in the same manner. We
show  again  just how to bound one of them:
\begin{eqnarray}
|T_5|\le
\|F\|_{L^2}\|\mathcal{K}(F-G)\|_{L^\infty}\|\partial_\theta(F-G)\|_{L^2}\nonumber\\
\le R_2\|\mathcal{K}\|_{L^2\to
L^\infty}\|F-G\|_{L^2}\|\partial_\theta(F-G)\|_{L^2}\le R_2
b\|k\|_{L^\infty}\sqrt{2\pi}\|\partial_\theta(F-G)\|_{L^2}^2\nonumber,
\end{eqnarray} where we used Poincar\'{e}'s inequality in the last
relation with the constant $R_2=\epsilon+(2s+b\|\partial_{\theta\theta}k\|_{L^\infty})(8\bar
C+\frac{6}{\pi})$, for some $\epsilon>0$  and the relation holds for the time
$t\ge t_1$ after which $F=\partial_\theta f\in B_{L^2}(0,R_2)$ (see
{\it Theorem}~\ref{dissipativity}).

\par Hence

\begin{equation}
|T_5+T_8|\le (R_2 +R_1)b\|k\|_{L^\infty}\sqrt{2\pi}
\|\partial_\theta(F-G)\|_{L^2}^2, \label{58}
\end{equation} where the inequality holds for $t\ge
t_2=\max\{t_1,t_0\}$, where $t_2$ is the time after which $f\in
B(0,R_1),F=\partial_\theta f\in B(0,R_2)$ (see {\it
Theorem}~\ref{dissipativity}).
\par Using bounds (\ref{1346}),(\ref{27}),(\ref{58}) into
(\ref{big}) together with Poincar\'{e}'s inequality and assumption
(\ref{condh1}) we obtain that $\|F-G\|_{L^2}^2$ decays exponentially
after time $t_2$.

\medskip\par{\it Step 2} Consider the ball $B(0,R)$ in $H^1$
for some  $R>\max \{R_1,R_2\}$. Then {\it
 Theorem} ~\ref{dissipativity} and the first
step show that there exists a time $t_3$ such that for all $t>t_3$
we have

\begin{equation}
S(t):B(0,R)\cap\mathcal{C}\cap S_{L^1}(0,1)\to
B(0,R)\cap\mathcal{C}\cap S_{L^1}(0,1), \label{ballinvariance}
\end{equation}
and there exists some $\alpha<1$, such that
\begin{equation}
\|S(t)f_0-S(t)g_0\|_{H^1}\le \alpha \|f_0-g_0\|_{H^1},
\label{contraction}
\end{equation} where $S(t)$ denotes the
nonlinear semigroup generated by the equation. We denoted by
$\mathcal{C}$ the cone of nonnegative functions and by
$S_{L^1}(0,1)$ the sphere of radius $1$, centered at $0$, in $L^1$.
\par Let $X=B(0,R)\cap\mathcal{C}\cap S_{L^1}(0,1)$. Then
$X$ with the metric induced by the $H^1$ norm is a complete metric
space. Define $T:X\to X,\,\, Tf=S(a)f$ for some $a\in\mathbb{Q},a>t_3$. The previous
arguments show that $T$ thus defined is a contraction.
\par Denote $T\circ T\circ\dots\circ T=T^n$. As $T$ is a
contraction we have that as $n\to\infty$, $T^nf_0\to f_1$ where
$f_1$ is a fixed point of $T$.
\par Similarly, taking $Uf=S(b)f$ for $b\in\mathbb{R}-\mathbb{Q},b>t_3$, we have
 $U:X\to X$ is a contraction. Reasoning as before we obtain the
 existence of a $f_2\in X$ such that $Uf_2=f_2$.
 \par From {\it Step 1} we have
 $$\lim_{n\to\infty}\|S(na)f_1-S(na)f_2\|_{H^1}=0.$$

 \par But $S(na)f_1=T^nf_1=f_1$ so the last limit becomes:
 \begin{equation}
 \|f_1-S(na)f_2\|_{H^1}\to 0,\,\textrm{as}\,\,n\to\infty.\label{eq:conta}
 \end{equation}

 \par Recall  Hurwitz's theorem in number theory (see
 for instance
 \cite{RADE64}) which states that  for $\gamma\in \mathbb{R}-\mathbb{Q}$ there
 are infinitely many rationals $\frac{p}{q}$ such that

 \begin{displaymath}
 |\gamma-\frac{p}{q}|<\frac{1}{\sqrt{5}q^2}.
 \end{displaymath}

 \par An easy consequence is that there exist two sequences
 $(m_k)_{k\in\mathbb{N}},(n_k)_{k\in\mathbb{N}};m_k,n_k\ge k,\forall
 k\in\mathbb{N}$, such that

 \begin{displaymath}
 |n_k a-m_k b|<\frac{1}{k}.
 \end{displaymath}

\par Let $\epsilon>0$. As $S(t)f_2$ is continuous
at $t=0$ there exists a $\delta>0$ such that

\begin{equation}
|t|<\delta\Rightarrow \|S(t)f_2-f_2\|_{H^1}<\frac{\epsilon}{2}.
\label{eq:contb}
\end{equation}

\par Thus, for $k$ large enough so that
$|n_ka-m_kb|<\frac{1}{k}<\delta$ we have

\begin{eqnarray}
\|S(n_k
a)f_2-S(m_kb)f_2\|_{H^1}=\|S(n_ka-m_kb)S(m_kb)f_2-S(m_kb)f_2\|_{H^1}
\nonumber\\
=\|S(n_ka-m_kb)f_2-f_2\|_{H^1}< \frac{\epsilon}{2}, \label{half1}
\end{eqnarray} where for the first equality we used the semigroup
property and for the last inequality relation (\ref{eq:contb}).
\par On the other hand from (\ref{eq:conta}) we know that for there
exists a rank $n_0$ such that for $n_k\ge n_0$

\begin{equation}
\|f_1-S(n_ka)f_2\|_{H^1}\le \frac{\epsilon}{2} \label{half2}.
\end{equation}

\par Putting together (\ref{half1}) and (\ref{half2}) we
obtain

\begin{displaymath}
\|f_1-f_2\|_{H^1}<\epsilon,
\end{displaymath} and since $\epsilon$ is arbitrary we have
that $f_1=f_2$.

\par We show now that $S(t)f_1$ is periodic of arbitrarily small
period.
\par Take $m_k,n_k$ such that $|n_ka-m_kb|<\frac{1}{k}$. Assuming
without loss of generality that $n_ka<m_kb$ we have
$$S(m_kb-n_ka)f_1=S(m_kb-n_ka)S(n_ka)f_1=S(m_kb)f_1=S(m_kb)f_2=f_2=f_1,$$
and thus $S(t)f_1$ has time period $0<m_kb-n_ka<\frac{1}{k}$. It is
well known that a continuous function of arbitrarily small periods
must be constant. Thus $f_1$ is a steady state.
\par As the difference of any two solutions tends to $0$ as
$t\to\infty$ this shows that all solutions tend to the steady state
$f_1$, which must be unique. \end{proof}
Following the proof, one can easily see that we also have

\begin{corollary} If (\ref{condl2}) and (\ref{condh1}) are satisfied
there exists a unique steady state solution of the equation
(\ref{eq:fred}), for an arbitrary smooth potential $k$ satisfying
(\ref{k}).
\end{corollary}

\smallskip
\par Let us observe now that the presence of the flow introduces in the equation
 a term of the form $\partial_\theta[(\omega+s\cos(2\theta)\bar
 f)]$.  The parameters $\omega$ and $s$ play very different r\^{o}les.
 In the case when either $\omega$ or $s$ is zero
 we have that the equation still has a gradient structure, in
 an appropriate reference frame.

\begin{lemma} 
\label{lemhomo1}Consider equation (\ref{eq:fred}).  The following holds:
\par {(i)} If $\omega=0,s\ne 0$  equation (\ref{eq:fred}) is of gradient type and
the $\omega$-limit set of any solution consists of steady states
solutions of (\ref{eq:fred}).
 \par{(ii)} If $s=0,\omega\ne 0$ make the rotating frame transformation
\begin{displaymath}
\tilde f(t,\theta)=\bar f(t,\theta+\omega t).
\end{displaymath}
\par Then $\tilde f$ satisfies the equation (\ref{eq:fgen}) (with
$V=0$) which  is an equation of gradient type.
\par Moreover, we have a time periodic solution (in  the moving frame) of (\ref{eq:fred}), namely
$g(\theta-\omega t)$ where $g$ is a steady state solution of
(\ref{eq:fred}) for $\omega=s=0$.
\smallskip
\par{(iii)} ({\rm an isotropic-nematic pattern}) Assume that $s=0,\omega\ne 0$.
 Consider solutions of (\ref{eq:fred}) with Maier-Saupe potential, for which
the initial data is even around $0$, in the $\theta$ variable. As
$t\to\infty$ we  have

\begin{displaymath}
\bar f(t,x_1,x_2,\theta)\stackrel{H^1}{\to g}(\theta-\omega t),
\end{displaymath} where $g=\frac{1}{2\pi}$ if $y_1(f(0,\cdot))=0$ and
$g\in\{f_{r(b)},f_{-r(b)}\}$  if $y_1(\bar f(0,\cdot))\ne 0$ (with
$r(b)$ sastisfies \eqref{Bessel2} and
$y_1(\bar f(0,\cdot))=\int_0^{2\pi} \bar f(0,x_1,x_2,\theta)\cos(2\theta)d\theta$).
\label{lemma:extreme}
\end{lemma}
\begin{proof} {(i)} In this case  condition (\ref{gradreq}) is satisfied,
as
$V(\theta)=\cos(2\theta)=\frac{d}{d\theta}\frac{\sin(2\theta)}{2}$.
Lemma~\ref{lemma:decaytosteady} gives us the conclusion.

\smallskip
\par{\it (ii)} Observe that $\mathcal{K}$ is invariant under
 the rotating frame transformation

\begin{equation}
(\mathcal{K}\bar f)(\theta+\omega t)=\int_{S^1}k(\theta+\omega
t-\theta')\bar f(\theta')d\theta' =\int_{S^1}
k(\theta+\omega t-\theta'-\omega t)\bar f(\theta'+\omega
t)d\theta'=(\mathcal{K}\tilde f)(\theta),\label{invariance}
\end{equation} thus $\tilde f$  satisfies
  the equation
 \begin{equation}
 \partial_t\tilde f+\partial_\theta(\tilde f\partial_\theta(\mathcal{K}\tilde f))=
 \partial_{\theta\theta}\tilde f,\label{eq:ftilde}
 \end{equation} which is of gradient type.
 \par A direct computation, using (\ref{invariance}), shows that
 $g(\theta-\omega t)$ is a solution of (\ref{eq:fred}), when $s=0$
 (where $g$ is a steady state solution of (\ref{eq:ftilde})).

\smallskip \par {\it (iii)}  Consider the same rotating frame
transformation as in the previous part. Using the fact that the
initial data is the same in the moving frame as in the fixed frame
and taking into account Theorem \ref{Thm2} we are
done. \end{proof}

\subsection{Maier-Saupe potential case $\omega=0, s\ne 0$: Stationary solutions and asymptotic behaviour}
Lemma \ref{lemhomo1} (i) has demonstrated that when $\omega=0$ and $s\ne 0$, a stationary solution of equation \eqref{eq:fred} always exists. In this subsection, we investigate deeper into the characteristics of the stationary solution, particularly focusing on the case of the Maier-Saupe potential and its large-time asymptotic behavior approaching this stationary state. It is worth recalling that the equation takes the following form:
\begin{equation}\label{equations1}
\partial_t \bar f+\partial_\theta[(s\cos(2\theta))\bar f+
\partial_\theta\mathcal{K}\bar f\cdot \bar f]=\partial_{\theta\theta}\bar f.
\end{equation}
Similarly as Section \ref{sec2}, the key tools are the free energy and Fisher information, defined as
\begin{equation}
\label{energys2}
\mathcal E(\bar f):=\int_0^{2\pi}\bar f\log\bar fd\theta-\frac b4\left(\int_0^{2\pi}{\bar f\sin(2\theta)d\theta}\right)^2-\frac b4\left(\int_0^{2\pi}{\bar f\cos(2\theta)d\theta}\right)^2-\frac s2\int_0^{2\pi}\bar f\sin(2\theta)d\theta,
\end{equation}
\begin{equation}
\label{fishers2}
\mathcal I[\bar f]:=-\frac d{dt}\mathcal E(\bar f)=\int_0^{2\pi}|\partial_\theta (\log\bar f-\mathcal K\bar f-\frac s2 \sin(2\theta))|^2f d\theta.
\end{equation}
We first give the form of the stationary solution  $\bar f$. Observe that $\bar f$ satisfies 
\begin{equation}\label{barf}
\partial_\theta(\mathcal{K}\bar f+\frac s2\sin(2\theta))\cdot \bar f=\partial_{\theta}\bar f,
\end{equation}
hence $\bar f$ has the form $\frac{e^{r\cos(2\theta)+q\sin(2\theta)}}{\int_0^{2\pi}e^{r\cos(2\theta)+q\sin(2\theta)}d\theta}$, where $r,q\in\mathbb R$ depend on $b,s$. Note that
$\int_0^{2\pi}e^{r\cos(2\theta)+q\sin(2\theta)}d\theta= 2\pi I_0(\sqrt{r^2+q^2})$. Denote $\varphi \in [0, 2\pi]$ that satisfies $\cos\varphi= \frac r{\sqrt{r^2+q^2}}, \sin\varphi= \frac q{\sqrt{r^2+q^2}}$. Then
\begin{equation*}
\begin{aligned}
&\int_0^{2\pi}e^{r\cos(2\theta)+q\sin(2\theta)}\cos(2\theta)d\theta=\int_0^{2\pi} e^{\sqrt{r^2+q^2}\cos(2\theta-\varphi)}\cos(2\theta) d\theta=\int_0^{2\pi} e^{\sqrt{r^2+q^2}\cos(2\theta)}\cos(2\theta+\varphi) d\theta\\
&=\frac r{\sqrt{r^2+q^2}}\int_0^{2\pi} e^{\sqrt{r^2+q^2}\cos(2\theta)}\cos(2\theta) d\theta-\underbrace{\frac q{\sqrt{r^2+q^2}}\int_0^{2\pi} e^{\sqrt{r^2+q^2}\cos(2\theta)}\sin(2\theta) d\theta}_{=0}\\
&= \frac{2\pi rI_1(\sqrt{r^2+q^2})}{\sqrt{r^2+q^2}}.
\end{aligned}
\end{equation*}
Similarly, \[ \int_0^{2\pi}e^{r\cos(2\theta)+q\sin(2\theta)}\sin(2\theta)d\theta=  \frac{2\pi qI_1(\sqrt{r^2+q^2})}{\sqrt{r^2+q^2}}.
\]
So we obtain that 
\[
\mathcal K\bar f=\frac b2 \frac{r}{\sqrt{r^2+q^2}}\frac{I_1(\sqrt{r^2+q^2})}{I_0(\sqrt{r^2+q^2})} \cos (2\theta)+\frac b2 \frac{q}{\sqrt{r^2+q^2}}\frac{I_1(\sqrt{r^2+q^2})}{I_0(\sqrt{r^2+q^2})} \sin (2\theta).
\]
From \eqref{barf}, $\log \bar f-\mathcal K\bar f-\frac s2\sin (2\theta)$ is constant. This means that for any $\theta\in [0, 2\pi]$, $r\cos(2\theta)+(q-\frac s2)\sin(2\theta)= \mathcal K\bar f$. Comparing the coefficients of $\cos(2\theta), \sin(2\theta)$, we obtain that $r,q$ satisfy 
\begin{equation}
\label{equationrq1}
\frac b2\frac{r}{\sqrt{r^2+q^2}}\frac{I_1(\sqrt{r^2+q^2})}{I_0(\sqrt{r^2+q^2})}=r, \quad
\frac b2\frac{q}{\sqrt{r^2+q^2}}\frac{I_1(\sqrt{r^2+q^2})}{I_0(\sqrt{r^2+q^2})}-q=-\frac s2.
\end{equation}
If $r\ne 0$, then from \eqref{equationrq1} we have $s=0$, a contradiction. So $r=0$. Thus \eqref{equationrq1} degenerates to the following equation 
\begin{equation}
\label{equationrq2}
q-\frac b2\frac{I_1(q)}{I_0(q)}=\frac s2.
\end{equation}
Denote $F(q):=q-\frac b2\frac{I_1(q)}{I_0(q)}$. Note that $F$ is an odd function of $q$, we can suppose that $s>0$. So finally we obtain that $\bar f$ has the form $\bar f_q(\theta):=\frac{e^{q\sin(2\theta)}}{\int_0^{2\pi}e^{q\sin(2\theta)}d\theta}$, where $q$ satisfies \eqref{equationrq2}. Moveover, we have the following lemma:
\begin{lemma}\label{bes} For any $q\in\mathbb R$, $(\frac{I_1(q)}{I_0(q)})'\in (0,\frac 12]$.
\end{lemma}
\begin{proof}
Notice that $I_0$ is even function, $I_1$ is odd, so $\left(\frac{I_1(z)}{I_0(z)}\right)'$ is even. From \eqref{Bessel}, we have $\frac{I_1(z)}{I_0(z)}\sim \frac z2$ as $z\to 0$, which means $\left(\frac{I_1(z)}{I_0(z)}\right)'|_{z=0}=\frac 12$.
 Thus we only need t prove $\left(\frac{I_1(z)}{I_0(z)}\right)'\in (0,\frac 12)$ for all $z>0$. It is well known that for any $z\in \mathbb R$, $\nu\in \mathbb Z$ (see (B.31), (B.41) of \cite{MAI} for example) ,
\[
I_{-\nu}(z)=I_\nu(z), \quad I_\nu(z)=\frac z{2\nu}(I_{\nu-1}(z)-I_{\nu+1}(z)), \quad I'_\nu(z)=\frac 12(I_{\nu-1}(z)+I_{\nu+1}(z)).
\]
So
\begin{equation*}
\begin{aligned}
\left(\frac{I_1(z)}{I_0(z)}\right)'=\frac{I_0(z)(I_0(z)+I_2(z))-I_1(z)(I_{-1}(z)+I_1(z))}{2I^2_0(z)}=\frac 1{I_0^2}\left( I_0^2-\frac1z I_0I_1-I_1^2\right).
\end{aligned}
\end{equation*}
Thus we need to show that
\[
I_0^2-\frac1z I_0I_1-I_1^2>0, \quad I_0^2-\frac2z I_0I_1-2I_1^2<0.
\]
the first inequality is because of $\xi(b)>0$ in \eqref{xi1} and $b=\frac{2rI_0(r)}{I_1(r)}$,  and the  second inequality is because of \eqref{xi2}.
\end{proof} From Lemma \ref{bes}, we obtain that when $0<b\le 4$, $F$ is an  increasing function of $q$. When $b>4$, there exists a unique $q_1(b)>0$, such that $F$ is decreasing  on $[0,q_1]$ and increasing on $[q_1,\infty)$, and moreover $F(q_1)<0$. Thus we have the following proposition about stationary solutions of  \eqref{eq:fred}.

\begin{proposition}
\label{propstationary}
 Consider equation  \eqref{equations1}.
 \\
(1) If $0<b\le 4$, then for any $s>0$, there exists a unique $q(b,s)>0$, such that $\bar f_q$ is the stationary solution of \eqref{equations1}.\\
(2) If $b>4$, suppose that $q_1>0$ satisfies $(\frac{I_1(q)}{I_0(q)})'=\frac 2b$, and $q_2>0$ is the positive zero of $F$, $q_3>0$ satisfies $F(q_3)=-F(q_1)$.\\
(2.1)For $s>-F(q_1)$,  \eqref{equations1} has one stationary solution $\bar f_q$, where $q>q_3$.\\
(2.2)For $s=-F(q_1)$,  \eqref{equations1} has two stationary solutions $\bar f_{-q_1}$ and $\bar f_{q_3}$.\\
(2.3)For $0<s<-F(q_1)$, \eqref{equations1} has three stationary solutions $\bar f_q, \bar f_{q'}, \bar f_{q''}$, where $q''\in (-q_2,-q_1), q'\in (-q_1,0), q\in (q_2,q_3)$.\\
\end{proposition}
Proposition \ref{propstationary} implies that the stationary solution of  \eqref{equations1} is not unique if $s\le -F(q_1)$. However, by comparing the free energy, we can show that
 $\bar f$ will converge to a specific stationary solution if  $\mathcal E[\bar f_{\text int}]$ is small enough.

\begin{proposition}
For $0<s<-F(q_1)$, $\mathcal E(\bar f_q)< \mathcal E(\bar f_{q''})<\mathcal E(\bar f_q')$, and  for $s=-F(q_1)$, $\mathcal E(\bar f_{q_3})<\mathcal E(\bar f_{-q_1})$.
\end{proposition}
\begin{proof}
We only prove the case $0<s<-F(q_1)$, the case $s=-F(q_1)$ is similar. For any $\alpha>0$, from direct calculation we have:
\[
\mathcal E(\bar f_\alpha)=-\log I_0(\alpha)+\alpha\frac{I_1(\alpha)}{I_0(\alpha)}-\frac b4\frac{I^2_1(\alpha)}{I^2_0(\alpha)}-\frac s2\frac{I_1(\alpha)}{I_0(\alpha)}.
\]
Recalling that for $\alpha=q,q',q''$, $\frac{I_1(\alpha)}{I_0(\alpha)}=\frac{2\alpha-s}{b}$. So we need to consider the function
\[
G(\alpha):=-\log I_0(\alpha)+\frac 1b(\alpha^2-s\alpha).
\]
Observing that $G'(\alpha)=-\frac {I_1(\alpha)}{I_0(\alpha)}+\frac 1b(2\alpha-s)=\frac 2b(F(\alpha)-\frac s2)$, we  deduce that $G(q')>G(q''), G(q')>G(q)$, and $G(q')-G(q)>G(q')-G(q'')$.
 This means that $G(q)<G(q'')<G(q')$.
\end{proof}

\begin{remark}
For $\bar f_q$ as the minimizer of $\mathcal E$, we could also study the convergence and asymptotic behaviour of $\bar f$ towards $\bar f_q$, and we could start focusing on the case that $\bar f=\bar f(\sin 2\theta)$. We could also similarly define quadratic forms $Q_1$ and $Q_2$ around $\bar f_q$ and the positiveness of $Q_1$ and the coercivity result between $Q_2$ and $Q_1$. Finally, we define the scalar product $\langle \cdot, \cdot\rangle$ as before, and we could show that  $\langle \bar g,\bar g\rangle=Q_1[\bar g], \langle \mathcal L\bar g,\bar g\rangle=-Q_2[\bar g]$. Moreover, the nonlinear term has higher order and could be controlled by $Q_1$ and $Q_2$, and we have the asymptotic behaviour by using Gr{\"o}nwall inequality. We skip the details here for simplicity.
\end{remark}

\subsection{Maier-Saupe ponential case $\omega, s\ne 0$: Existence of time periodic solution}
\par Lemma \ref{lemhomo1} suggests that one can think, heuristically, that
 $\omega$ affects the rotational
behavior of the solution while $s$ affects the long time shape of
the solution. The r\^{o}les of $\omega$ and $s$ are thus very
different. When they are both non-zero one can not reduce the
equation to a gradient type one just by considering it in the
rotating frame. What prevents us from repeating the argument  above
is the fact that in a moving frame the linear potential $V$  becomes
time dependent. In this situation there is no obvious Lyapunov
functional.

\par Nevertheless, when both $s$ and $\omega$
are nonzero, we will be able to show that for small $s$ a  time-periodic solution of a period $\frac{\pi}{\omega-s}$ does not exist in a neighbourhood of the periodic one existing at $s=0$. However this leave open the possibility that there still exist time periodic solutions with different periods or further away from the known periodic solution. 
\par The following argument is done just for the Maier-Saupe
potential. This is because at the present time only for this
potential there is a good understanding of the form of the steady
state solutions (in the absence of the flow) and of their dependence
on the concentration intensity parameter $b$.

\bigskip

\begin{theorem} \label{thm:periodic} Consider  equation (\ref{eq:fred}) with $\omega>0
$ and let $\mathcal{K}=\mathcal{K}_{MS}$ be the Maier-Saupe
potential. For $b>4$   and $s\not=0$ small enough (depending on $\omega$) the equation (\ref{eq:fred}) has no-time periodic solution of period $\frac{\pi}{\omega-s}$ in the neighbourhood of the time periodic solution corresponding to formally setting $s=0$ (as provided in  Lemma~\eqref{lemhomo1}). \label{theorem:periodic}
\end{theorem}

\begin{proof} We  present first the strategy of the proof. We are interested in
obtaining zeroes of a
 functional $\mathcal{F}:\mathcal{X}\times \mathbb{R}\times\mathbb{R}\to\mathcal{Y}$

\begin{displaymath}
\mathcal{F}(f,\omega,s)=\partial_t
f+\partial_\theta\left[(\omega+s\cos(2\theta))f\right]+\partial_\theta\left[
\partial_\theta(\mathcal{K}f)\cdot f\right]-\partial_{\theta\theta}f,
\end{displaymath} where $\mathcal{X},\mathcal{Y}$ are spaces of
functions periodic both in $t$ and $\theta$, whose precise
definition will not be given because we will see soon that it is
more convenient to work with a different formulation of the above
functional. For that formulation we will make precise the functional
spaces.
\par We will show that for any $\omega$ there exists a
$\lambda_{\omega}$ such that for any $\lambda\in
(-\lambda_{\omega},\lambda_{\omega})$ we have

\begin{equation}
\mathcal{F}(f,\omega+\lambda,\lambda)=0\label{transition}
\end{equation} for some $f\in\mathcal{X}$. This
 suffices for obtaining the conclusion of the theorem. 
 
Returning to (\ref{transition}) let us consider
the  Smoluchowski equation in homogeneous flow
\begin{equation}
\partial_t
f+\partial_\theta[(\omega+\lambda
+\lambda\cos(2\theta))f+\partial_\theta(\mathcal{K}f)f]=\partial_{\theta\theta}
f. \label{eq:SmolS}
\end{equation}
Make the rotating frame transformation $\tilde
f(t,\theta)=f(t,\theta+\omega t)$. Then (\ref{eq:SmolS}) becomes
\begin{equation}
\partial_t \tilde f+\lambda\partial_\theta [\tilde f+\cos(2\theta+2\omega t)\tilde
f]+\partial_\theta [\partial_\theta(\mathcal{K}\tilde f)\cdot \tilde
f]=\partial_{\theta\theta} \tilde f. \label{eq:SmolSmoving}
\end{equation}
Let $g(\theta)$ be an even,  {\it nonconstant} solution of $\partial_\theta[\partial_\theta(\mathcal{K}f)f]=\partial_{\theta\theta}
f$. We know that such a solution exists and it is given by the
formula $g=\frac{e^{r(b)\cos(2\theta)}}{I_0(r(b))}$ for a certain $r(b)$ satisfying \eqref{Bessel2}. Next, we  decompose $\tilde f(t,\theta)=z(t,\theta)+g$. Then $z$ satisfies the equation
\[\partial_t z+\lambda\partial_\theta [(z+g)+\cos(2\theta+2\omega t)(z+g)]+
\partial_\theta[\partial_\theta(\mathcal{K}z)g+
\partial_\theta(\mathcal{K}g) z+\partial_\theta(\mathcal{K}z) z]=z_{\theta\theta}.\]
We  prove the existence of time periodic solutions for the above
equation, with time period $\frac{\pi}{\omega}$. Define
$F:\mathcal{X}\times \mathbb{R}\to \mathcal{Y}$ by
\begin{equation}
F(z,\lambda)=\partial_t z+\lambda\partial_\theta
[(z+g)+\cos(2\theta+2\omega t)(z+g)]
+\partial_\theta[\partial_\theta(\mathcal{K}z)g+
\partial_\theta(\mathcal{K}g) z+\partial_\theta(\mathcal{K}z) z]-z_{\theta\theta}
\label{def:F}
\end{equation} with
\begin{eqnarray}
\mathcal{X}=\{z\in
H^1\left([0,\frac{\pi}{\omega}],H^2[0,2\pi]\right), z(t,0)=z(t,2\pi),\forall t\in [0,\frac{\pi}{\omega}]
,\nonumber\\z(0,\theta)=z(\frac{\pi}{\omega},\theta),\partial_t z(0,\theta)=\partial_t z(\frac{\pi}{\omega},\theta,)
\forall \theta\in
[0,2\pi],
\nonumber\int_0^{2\pi}z(t,\theta)d\theta=0,\forall
t\in [0,\frac{\pi}{\omega}]\}\nonumber
\end{eqnarray} and
\begin{eqnarray}
\mathcal{Y}=\{z\in
L^2\left([0,\frac{\pi}{\omega}],L^2[0,2\pi]\right),
z(0,\theta)=z(\frac{\pi}{\omega},\theta),\forall \theta\in
[0,2\pi],\nonumber\\
 z(t,0)=z(t,2\pi),\forall t\in [0,\frac{\pi}{\omega}],\int_0^{2\pi}z(t,\theta)d\theta=0,\forall
t\in [0,\frac{\pi}{\omega}]\}\nonumber.
\end{eqnarray}
\begin{remark} The choice of regularity spaces is somewhat arbitrary, as one
can see a posteriori that the solution is analytic in time and
space. What we need for our proof is that  the norm of $\mathcal{X}$
controls the $L^\infty$ norm in time and space.  We also need that
the structure of $\mathcal{Y}$ allows for a simple orthogonal
decomposition.
\end{remark}

\par We have then that

\begin{displaymath}
F(0,0)=0.
\end{displaymath}
We want to apply the implicit function theorem and obtain the
existence of a periodic solution for small $\lambda$. This is a
continuation argument, which finds a periodic solution of period
$\frac{\pi}{\omega}$ near one which we already know to exists (for
$\lambda=0$, see {\it Lemma}~\ref{lemma:extreme}). It is due to the
fact that we are working in a rotating frame that the time periodic
 solution in the initial frame, $g(\theta-\omega t)$, is stationary in
  the rotating frame, $g$ in (\ref{def:F}).

In order to apply the implicit function theorem we need thus to check
that $$Lh=\partial_z F(0,0)h=\partial_t
h+\partial_\theta[\partial_\theta\mathcal{K}g\cdot
h+\partial_\theta\mathcal{K}h\cdot g]-\partial_{\theta\theta}h$$ as
a bounded operator from $\mathcal{X}$ to $\mathcal{Y}$ is a
homeomorphism, i.e., taking into account the open mapping theorem,
that $L$ is bijective. Nevertheless, this is not the case since, as
we will see, $\textrm{dim}(ker(L))=\textrm{codim}(range(L))=1$, so $L$ is a an operator
of Fredholm index $0$. In this situation  an implicit function
theorem is still possible under a certain {\it "non-resonance"}
condition. This is available for instance in \cite{HIKI04}, p.12. We
will present it in {\it Lemma~\ref{fredh0}} below after analyzing
the operator $L$.

\par In order to determine the kernel and the range of the operator
$L$ we need to study equations of the form $Lh=f$, for
$f\in\mathcal{Y}$. Multiplying such an equation by
$\sqrt{\frac{\omega}{\pi}}e^{-ik2\omega t}$ and integrating on
$[0,\frac{\pi}{\omega}]$ we have that the equation $Lh=f$ reduces to
a decoupled system of ordinary differential equations:

$$i2\omega k h_k+\partial_\theta
[\partial_\theta(\mathcal{K}g)h_k+\partial_\theta(\mathcal{K}h_k)g]-
\partial_{\theta\theta}h_k=f_k,$$ where we denote by $h_k,f_k$ the
$k$-th Fourier mode in time of $h$, respectively $f$.

\par Thus the problem of determining the kernel and range of $L$ reduces  to
understanding  the  operator
\begin{equation}\tilde L(h)=\partial_{\theta\theta}h-\partial_\theta
[\partial_\theta(\mathcal{K}g)h+\partial_\theta(\mathcal{K}h)g].
\label{tildeL}
\end{equation}
We have
\begin{lemma} Let $\tilde L:\tilde X\to \tilde Y$ be a bounded
operator as defined in (\ref{tildeL}) with $$\tilde X=\{f\in
H^2[0,2\pi],\int_0^{2\pi}f(\theta)d\theta=0\},\,\,\tilde Y=\{f\in
L^2[0,2\pi],\int_0^{2\pi}f(\theta)d\theta=0\}.$$
\par Then
\begin{displaymath}
\ker(\tilde L)= \{p\partial_\theta g,p\in\mathbb{R}\},
\end{displaymath} and
\begin{equation}
f\in Range(\tilde L)\Leftrightarrow \int_0^{2\pi}(\int_0^\theta
f(\sigma)d\sigma)(\frac{1}{g(\theta)}-\frac{1}{2\pi}\int_0^{2\pi}
\frac{d\sigma}{g(\sigma)})d\theta=0. \label{rangetildeL}
\end{equation}

\par Thus $\dim(Ker(\tilde L))=\operatorname{codim}(Range(\tilde L))=1$, so $\tilde
L$ is an  operator of Fredholm index $0$.
\par Moreover, regarding  $\tilde L$ as an unbounded operator on $\tilde Y$,
 with $D(\tilde L)=\tilde X$, we have that
$\tilde L$  is a sectorial operator
 with discrete spectrum   contained in the real line.
\end{lemma}

\par \begin{proof} Standard arguments show that the operator, regarded as an unbounded operators on $\tilde Y$,
 has discrete spectrum and it is sectorial (see for instance \cite{HENRY81},\cite{KATO80}). The proof will not be given here.
 Let us observe that the operator can only have real
 spectrum, i.e. only real eigenvalues. Define
\begin{equation}
Ah=\frac{h}{g}-Kh, \label{def:Ah}
\end{equation} and then

\begin{equation}
\tilde Lh=(h_\theta-h\frac{\partial_\theta g}{g}-
g(Kh)_\theta)_\theta+(\frac{h}{g} \underbrace{[\partial_\theta
g-g(Kg)_\theta]}_{=0})_\theta =(g(Ah)_\theta)_\theta,\label{eq:tildeL}
\end{equation} where the cancellation is due to our choice of $g$.

\par In order to compute the spectrum of $\tilde L$ consider the
equation

\begin{displaymath}
\tilde L(h_R+ih_I)=(R+iI)(h_R+ih_I),
\end{displaymath} with $R,I,h_R,h_I$ real quantities.

\par Separating the real and imaginary parts in the above equation
we obtain

\begin{equation}
(g(Ah_R)_\theta)_\theta=Rh_R-Ih_I, \label{eq:real}
\end{equation}
and
\begin{equation}
(g(Ah_I)_\theta)_\theta=Ih_R+Rh_I. \label{eq:imaginary}
\end{equation}

\par Multiplying (\ref{eq:real}) by $Ah_R$, adding the result
to (\ref{eq:imaginary}) multiplied by $Ah_I$, integrating over
$[0,2\pi]$ and by parts we have

\begin{equation}
-\int_0^{2\pi} g((Ah_R)_\theta)^2-\int_0^{2\pi}
g((Ah_I)_\theta)^2=R(Ah_R,h_R)+ R(Ah_I,h_I), \label{eq:cancel1}
\end{equation} where we used the fact that
\begin{equation}
\int_0^{2\pi} h_I Ah_R=\int_0^{2\pi}h_R Ah_I.\label{eq:Aselfadj}
\end{equation}

\par Also, let us observe that by multiplying (\ref{eq:real}) by
$Ah_I$, integrating over $[0,2\pi]$ and by parts, we obtain on the
left hand side of the equality the same thing as
multiplying (\ref{eq:imaginary}) by $Ah_R$, integrating over $[0,2\pi]$
and by parts. This implies the equality of the corresponding right
hand sides, i.e.

\begin{equation}
\int_0^{2\pi} Rh_RAh_I-\int_0^{2\pi} I h_IAh_I= \int_0^{2\pi} Ih_R
Ah_R+\int_0^{2\pi}R h_I Ah_R. \label{eq:cancel2}
\end{equation} And using again (\ref{eq:Aselfadj}) we have

\begin{equation}
-I(Ah_I,h_I)=I(h_R,Ah_R) \label{eq:cancel3},
\end{equation} which, for $I\ne 0$, used in (\ref{eq:cancel1})
implies

\begin{displaymath}
-\int_0^{2\pi} g((Ah_R)_\theta)^2-\int_0^{2\pi} g((Ah_I)_\theta)^2=0.
\end{displaymath}
\par Taking into account that $g> 0$ the last equality implies that
$(Ah_I)_\theta=(Ah_R)_\theta\equiv 0$. Using this into
(\ref{eq:real}),(\ref{eq:imaginary}) we obtain $h_I=h_R=0$. Thus
necessarily $I=0$, so  the imaginary part of an eigenvalue must be
zero.

\smallskip\par In order to compute the range and the kernel, take
$f\in \tilde Y$ and denote $F(\theta)=\int_0^\theta
f(\sigma)d\sigma$. Then

$$\tilde Lh=\partial_\theta(g(Ah)_\theta)=\partial_\theta F$$
implies

\begin{equation}
g(Ah)_\theta=F(\theta)+c_1.\label{eq:rel0}
\end{equation}

\par In order to determine $c_1$ we divide by $g$ on both sides of
the last equality and integrate on $[0,2\pi]$ obtaining

\begin{equation}
c_1=-\frac{\int_0^{2\pi}
\frac{F(\theta)}{g(\theta)}d\theta}{\int_0^{2\pi}
\frac{1}{g(\theta)}d\theta}. \label{eq:rel1}
\end{equation}

\par Also, integrating on $[0,2\pi]$ both sides of (\ref{eq:rel0})
we get

\begin{equation}\int_0^{2\pi}
g(Ah)_\theta=\int_0^{2\pi}F(\theta)d\theta+c_1\cdot 2\pi.
\label{eq:rel01}
\end{equation}

\par On the other hand

\begin{equation*}
\begin{aligned}
\int_0^{2\pi} g(Ah)_\theta=-\int_0^{2\pi} g_\theta\cdot
Ah\nonumber
&=2r(b)\int_0^{2\pi} g(\theta)\sin(2\theta)[\frac{h}{g}-
\frac{b}{2}
c(h)\cos(2\theta)-\frac{b}{2}s(h)\sin(2\theta)]d\theta\nonumber\\&
=2r(b)[s(h)-\frac{b}{2}(\frac{1}{2}-\frac{1}{2}(1-\frac{4}{b}))s(h)]=0,
\nonumber
\end{aligned}
\end{equation*} where we used in the first equality an integration by parts
and in the second the relation $g_\theta=-2r(b)g\sin(2\theta)$(see
the definition of $g$). Also we denoted
$c(h)=\int_0^{2\pi}h(\theta)\cos(2\theta)d\theta$,
$s(h)=\int_0^{2\pi}h(\theta)\sin(2\theta)d\theta$ and  we used the
fact that $\int_0^{2\pi}g(\theta)\cos(4\theta)d\theta=1-\frac{4}{b}$
(see \cite {CKT04}). Using the last computation in (\ref{eq:rel01})
we obtain

\begin{equation}
2\pi c_1=-\int_0^{2\pi} F(\theta)d\theta \label{eq:rel2}.
\end{equation}
\par From (\ref{eq:rel1}) and (\ref{eq:rel2}) we have that $c_1$ is
equal to both sides of the equality

$$-\frac{\int_0^{2\pi}
\frac{F(\theta)}{g(\theta)}d\theta}{\int_0^{2\pi}
\frac{1}{g(\theta)}d\theta}=-\frac{\int_0^{2\pi}F(\theta)d\theta}{2\pi},
$$ which implies a restriction on $F$ namely (\ref{rangetildeL}).

\par Observe that  relation (\ref{rangetildeL}) is the
only restriction on the range of $\tilde L$. Indeed, returning to
(\ref{eq:rel0}) and dividing by $g$, integrating and recalling the
definition of $Ah$, (\ref{def:Ah}) we have

\begin{equation}
h(\theta)=\frac{b}{2}c(h)g(\theta)\cos(2\theta)
+\frac{b}{2}s(h)g(\theta)\sin(2\theta)+g(\theta)\int_0^{\theta}
\frac{F(\sigma)+c_1}{g(\sigma)}d\sigma+c_2g(\theta), \label{reph}
\end{equation} where $c_2$ is a constant of integration to be
determined.

\par Multiplying (\ref{reph}) by $\cos(2\theta)$ and integrating
over $[0,2\pi]$ we obtain

$$c(h)=\frac{b}{2}c(h)[\frac{1}{2}+\frac{1}{2}(1-\frac{4}{b})]+
\int_0^{2\pi} g(\theta)\cos(2\theta)\left(\int_0^\theta
\frac{F(\sigma)+c_1}{g(\sigma)}d\sigma\right)d\theta+ c_2 c(g),
$$ which implies

\begin{equation}
c(h)(2-\frac{b}{2})=\int_0^{2\pi}
g(\theta)\cos(2\theta)\left(\int_0^{\theta}
\frac{F(\sigma)+c_1}{g(\sigma)}d\sigma\right)d\theta+c_2 c(g). \label{ch1}
\end{equation}

\par Also, integrating (\ref{reph}) over $[0,2\pi]$ and using the
fact that we are looking for solutions $h$ which have mean zero we
get

\begin{equation}
-c_2=\frac{b}{2}c(h)c(g)+\int_0^{2\pi}
g(\theta)\left(\int_0^{\theta}\frac{F(\sigma)+c_1}{g(\sigma)}d\sigma\right)d\theta.
\label{c2}
\end{equation}

\par Multiplying the last relation by $-c(g)$ and replacing the
expression for $c_2c(g)$ thus obtained into (\ref{ch1}) we have

\begin{equation}
c(h)(2-\frac{b}{2}+\frac{b}{2}c(g)^2)=-c(g)\int_0^{2\pi} g(\theta)
\left(\int_0^\theta
\frac{F(\sigma)+c_1}{g(\sigma)}d\sigma\right)d\theta+
\int_0^{2\pi} g(\theta)\cos(2\theta)\left(\int_0^\theta
\frac{F(\sigma)+c_1}{g(\sigma)}d\sigma\right)d\theta. \label{ch}
\end{equation}

\par In the last relation we can divide by
$2-\frac{b}{2}+\frac{b}{2}c^2(g)$ (which is nonzero as
$c^2(g)>1-\frac{4}{b}$, for $b>4$ see \cite {LIZA205}, {\it Theorem}
$2.1$) and thus we obtain an expression for $c(h)$ only in terms of
$F$ and $g$ (as $c_1$ can also be determined in terms of $F$ and $g$,
see (\ref{eq:rel1})).

\smallskip
\par Let us check that we can take $s(h)$ to be arbitrary in the
representation formula for $h$. Indeed, multiplying (\ref{reph}) by
$\sin(2\theta)$ and integrating over $[0,2\pi]$ we get

$$s(h)=\frac{b}{2}s(h)\underbrace{[\frac{1}{2}-\frac{1}{2}(1-\frac{4}{b})]}_{2/b}+
\int_0^{2\pi} g(\theta)\sin(2\theta)(\int_0^\theta
\frac{F(\sigma)+c_1}{g(\sigma)}d\sigma)d\theta, $$ which is always
true as

\begin{eqnarray}\int_0^{2\pi}\underbrace{g(\theta)\sin(2\theta)}_{=-g_\theta/(2r(b))}
(\int_0^\theta \frac{F(\sigma)+c_1}{g(\sigma)}d\sigma)d\theta=
-\frac{1}{2r(b)}g(2\pi)\underbrace{\int_0^{2\pi}(\frac{F(\sigma)+c_1}{g(\sigma)})d\sigma}_{=0}
\nonumber
+\frac{1}{2r(b)}\underbrace{\int_0^{2\pi}g(\theta)\frac{F(\theta)+c_1}{g(\theta)}d\theta}_{=0}
\nonumber,
\end{eqnarray}
where we used an integration by parts for the equality and
(\ref{eq:rel1}),(\ref{eq:rel2}) for the cancelations.

\par Summarizing: for a given $f\in \tilde L$ satisfying the
compatibility condition (\ref{rangetildeL}) we define $c_1$ by
relation (\ref{eq:rel1}) and use this in (\ref{ch}) to obtain an
expression for $c(h)$. We use this to get $c_2$ from (\ref{c2}) and
plug everything in (\ref{reph}). The $h$ thus obtained will satisfy
the equation $\tilde Lh=f$ for an arbitrary $s(h)$.
\par In particular, if we look for a solution of $\tilde Lh =0$ we
obtain that $h\in \{p\partial_\theta g,p\in\mathbb{R}\}$ which gives
us the kernel of $\tilde L$.

\end{proof}

 \par The properties of $\tilde L$  imply
 that the full operator $L$ has  the same kernel as $\tilde L$(where now
 $p\partial_\theta g,p\in\mathbb{R}$ is regarded as an element in
 $\mathcal{X}$). Also,
 the compatibility condition for $f\in\mathcal{Y}$ to
 be in the range of $L$ is
\begin{equation}
\int_0^{\pi/(\omega)} \int_0^{2\pi}(\int_0^\theta
f(\sigma,s)d\sigma)(\frac{1}{g(\theta)}-\frac{1}{2\pi}\int_0^{2\pi}
\frac{d\sigma}{g(\sigma)})d\theta ds=0. \label{comp}
 \end{equation} So $range(L)$ has
codimension one.

\par A technical remark is necessary at this point: we have that for
$f\in\mathcal{Y}$ satisfying the appropriate compatibility condition
as above there exists some $h$ such that $Lh=f$. In order to make
sure that this $h$ is in $\mathcal{X}$ we need to check the
regularity in time. While our argument so far does not give but
$L^2$ regularity in time, higher regularity in time is nevertheless
available. This is a consequence of the parabolic nature of the
operator $L$.

 \par More precisely, the equation

$$ik2\omega h_k+\tilde Lh_k=f_k$$ and the fact that $ik2\omega$ is in the resolvent
of $\tilde L$, imply that $h_k=(ik2\omega+\tilde L)^{-1}f_k,\,\forall k\in\mathbb{Z}-\{0\}$, with

\begin{equation}
\|h_k\|_{H^2}\le \|(ik2\omega+\tilde L)^{-1}\|_{L^2\to
H^2}\|f_k\|_{L^2}. \label{resolvnonzero}
\end{equation}

\par If $k=0$ we have (by the compatibility condition (\ref{comp})) $f_0\in\, range(\tilde L) $ and thus $h_0=\tilde L^{-1} f_0$, with

\begin{equation}
\|h_0\|_{H^2}\le \|\tilde L^{-1}\|_{L^2\to H^2}\|f_0\|_{L^2}.
\label{resolvzero}
\end{equation}

\par As $\tilde L$ is sectorial we have that for $|k|$ large enough
$|k|>k_0>0$ the number  $ik2\omega$ is in the resolvent (which we
already knew from the lemma) and moreover (see \cite{FRIED69})

\begin{equation}
\|(ik2\omega+\tilde L)^{-1}\|_{L^2\to H^2}\le C, \,\forall |k|>k_0>0,
\label{klarge}
\end{equation} where $C$ is a constant independent of $k$.
\par Let $\tilde C=\max\{C,\|(ik2\omega+\tilde L)^{-1}\|_{L^2\to
H^2},k=0,\pm 1,\pm 2,\dots, \pm k_0\}$. We have then

$$\|h_k\|_{H^2}\le \tilde C\|f_k\|,\forall k\in\mathbb{Z},$$ which
implies that $h\in L^2((0,\frac{\pi}{\omega}),H^2[0,2\pi])$, i.e.
the same regularity in time as $f$. Also, using the equation
$Lh=\partial_t h+\partial_\theta[\partial_\theta\mathcal{K}g\cdot
h+\partial_\theta\mathcal{K}h\cdot g]-\partial_{\theta\theta}h=f$ we
get $\partial_t h\in L^2((0,\frac{\pi}{\omega}),L^2(0,2\pi))$ so
$h\in C([0,\frac{\pi}{\omega}],L^2(0,2\pi))$  (see for instance
\cite{EVANS02}).
\par  On the other hand   we also  have that $h$ is a weak
solution of a Cauchy problem for the equation $Lh=f$, with initial
data $h(0)\in L^2(0,2\pi)$. But $h(0)=h(\frac{\pi}{\omega})$ by time
periodicity  and then the parabolic regularization effect implies
$h(0)=h(\frac{\pi}{\omega})\in H^2[0,2\pi]$. For $f\in
L^2(0,\frac{2\pi}{\omega},L^2(0,2\pi))$ the Cauchy problem with
initial data $h(0)\in H^2[0,2\pi]$ has a unique solution $g\in
H^{1}([0,\frac{2\pi}{\omega}],H^2[0,2\pi])$ with $g(0)=h(0)$ (see
\cite{FRIED69}). Moreover it can be shown that the uniqueness holds
for solutions which are only in
$C([0,\frac{\pi}{\omega}],L^2(0,2\pi))$. Then $g\equiv h$ and thus
$h$ will have the necessary regularity in time for being in $\mathcal{X}$.

\par The abstract lemma ( which extends the Implicit
Function Theorem) that we  need is

\begin{lemma} (\cite{HIKI04}, p.12)
Let $F:U\times V\to Z$ with $U\subset X,V\subset \mathbb{R}$, where
$X$ and $Z$ are Banach spaces. Assume that $F\in C^1(U\times V,Z)$
and:
\par $\bullet$ $F(x_0,\lambda_0)=0$ for some $(x_0,\lambda_0)\in U\times
V$, $Range(D_xF(x_0,\lambda_0))$ is closed in $Z$  and
$$\dim(Ker(D_xF(x_0,\lambda_0)))=\operatorname{codim}(Range(D_xF(x_0,\lambda_0)))=1.$$
\par $\bullet$ We have the {\it "non-resonance"} condition

\begin{equation}
D_\lambda F(x_0,\lambda_0)\not\in Range(D_xF(x_0,\lambda_0)).
\label{nonresonance}
\end{equation}

\par Then there exists a continuously differentiable curve through
$(x_0,\lambda_0)$

$$\{(x(r),\lambda(r))|r\in
(-\delta,\delta),(x(0),\lambda(0))=(x_0,\lambda_0)\},$$ such that

$$F(x(r),\lambda(r))=0,\,\textrm{for}\,\,r\in(-\delta,\delta)$$ for
some $\delta>0$ and all the solutions of $F(x,\lambda)=0$ in a
neighborhood of $(x_0,\lambda_0)$ belong to the curve specified
above. \label{fredh0}
\end{lemma}

\par The {\it "non-resonance"} condition (\ref{nonresonance})
in our case is

\begin{displaymath}
\int_0^{2\pi}
g(\theta)(\frac{1}{g(\theta)}-\frac{1}{2\pi}\int_0^{2\pi}
g^{-1}(\sigma)d\sigma)d\theta\ne 0,
\end{displaymath} or

\begin{equation}
\label{nonre1}(2\pi)^2\ne \int_0^{2\pi} g(\sigma)d\sigma\cdot \int_0^{2\pi}
g^{-1}(\sigma)d\sigma.
\end{equation}
Recalling that $g(\theta)=\frac{e^{r(b)\cos(2\theta)}}{Z(b)}$
with $Z(b)=\int_0^{2\pi} e^{r(b)\cos(2\theta)}d\theta$. From Cauchy-Schwartz inequality, \eqref{nonre1} means that $g$ is not a constant function, which is equivalent to $r(b)\ne 0$. Remind the result of Section \ref{sec2}, this is always true for $b>4$.









 We note that what we obtained so far is that all the zeroes of $F$ in a neighbourhood of $(0,0)$ are located on a curve $(z(r),\lambda(r))$ for $r\in (-\delta,\delta)$. One can explicitly obtain what this curve is, namely $z(t,r)=g(r+\theta)-g(\theta)$ and $\lambda(r)=0$. Indeed, this can be checked by plugging this choice into the definition \eqref{def:F} of $F$ and taking into account that $g(r+\theta)$ is a solution of $\partial_\theta(f\partial_\theta\mathcal{K}f)=\partial_{\theta\theta} f$ for arbitrary $r$, thanks to the rotational invariance of the equation. 

What we obtained shows in particular that in a neighbourhood of the solution $(g(\theta+\omega t),0)$ there is no time periodic solution for small $s\not=0$.

\end{proof}

\par{\bf Acknowledgements:} Stimulating discussions of A.Z. with Professors P. Polacik and J.
Vukadinovic are gratefully acknowledged. X.Li thanks the hosting of BCAM when he visited.
 X.Li and A.Zarnescu have been partially  supported by the Basque Government through the BERC 2022-2025
program and by the Spanish State Research Agency through BCAM Severo Ochoa excellence accreditation
CEX2021-01142-S funded by\\ MICIU/AEI/10.13039/501100011033 and through Grant PID2023-146764NB-I00
funded by\\ MICIU/AEI/10.13039/501100011033 and cofunded by the European Union. X.Li has been partially supported by the LTC-Transmath project funded by Fundac\'ion Euskampus. A.Zarnescu has been  also partially supported by a grant of the Ministry of Research, Innovation and Digitization, CNCS-UEFISCDI, project number PN-IV-P2-2.1-T-TE-2023-1704, within PNCDI IV.  The authors would like to thank the anonymous referees for careful reading of the paper and especially for pointing out certain limitations in the initial version of Theorem~\ref{thm:periodic} .

\end{document}